\documentclass[pdflatex,sn-mathphys,Numbered]{sn-jnl}

\usepackage{graphicx}
\usepackage{multirow}
\usepackage{amsmath,amssymb,amsfonts}
\usepackage{amsthm}
\usepackage{mathrsfs}
\usepackage[title]{appendix}
\usepackage{xcolor}
\usepackage{textcomp}
\usepackage{manyfoot}
\usepackage{booktabs}
\usepackage{algorithm}
\usepackage{algorithmicx}
\usepackage{algpseudocode}
\usepackage{listings}

\numberwithin{equation}{section}

\newcommand{\beq}{\begin{equation}}
\newcommand{\eeq}{\end{equation}}
\newcommand{\beqa}{\begin{eqnarray}}
\newcommand{\eeqa}{\end{eqnarray}}

\newcommand{\be}{\begin{equation}}
\newcommand{\ee}{\end{equation}}
\newcommand{\bea}{\begin{eqnarray}}
\newcommand{\eea}{\end{eqnarray}}

\newcommand{\abs}[1]{\vert#1\vert}
\newcommand{\cc}{{\cal C}}
\newcommand{\dd}{{\rm d}}
\newcommand{\erf}{\mathop{\rm erf}}
\newcommand{\erfc}{\mathop{\rm erfc}}
\newcommand{\half}{{\scriptstyle{\frac{\scriptstyle1}{\scriptstyle2}}}}
\newcommand{\ii}{{\rm i}}
\newcommand{\mean}[1]{\langle#1\rangle}
\newcommand{\prob}{\mathbb{P}}

\newcommand{\lam}{{\lambda}}
\newcommand{\nn}{{\mathcal{N}}}
\newcommand{\nb}{{\mathcal{N}}_t^{\bullet}}
\newcommand{\nc}{{\mathcal{N}}_t^{\times}}
\newcommand{\ncminus}{{\mathcal{N}}_{t-1}^{\times}}
\newcommand{\ncb}{{\mathcal{N}}_t^{\times\bullet}}

\newcommand{\ib}{I^\bullet}
\newcommand{\ic}{I^{\times}}
\newcommand{\ab}{A^\bullet}
\newcommand{\ac}{A^\times}

\newcommand{\rw}{\breve{w}}
\newcommand{\stwo}[2]{\left\{#1\atop#2\right\}}

\newcommand{\A}{{\sf A}}
\newcommand{\B}{{\sf B}}
\newcommand{\C}{{\sf C}}
\newcommand{\Y}{{\cal Y}}
\newcommand{\Z}{{\cal Z}}

\newcommand{\T}{\boldsymbol T}
\newcommand{\mzeta}{{\boldsymbol\zeta}}
\newcommand{\ttau}{{\boldsymbol{\tau}}}
\newcommand{\Tz}{\boldsymbol{T}_{0\to0}}
\newcommand{\e}{{\rm e}}

\begin{document}

\title{Returns to the origin of the P\'olya walk with stochastic resetting}

\author[]{\fnm{Claude} \sur{Godr\`eche}}\email{claude.godreche@ipht.fr}

\author*[]{\fnm{Jean-Marc} \sur{Luck*}}\email{jean-marc.luck@ipht.fr}

\affil[]
{\orgdiv{Universit\'e Paris-Saclay, CEA, CNRS},
\orgname{Institut de Physique Th\'eorique},
\postcode{91191}
\city{Gif-sur-Yvette},
\country{France}}

\abstract{
We consider the simple random walk (or P\'olya walk) on the one-dimensional lattice
subject to stochastic resetting to the origin with probability $r$ at each time step.
The focus is on the joint statistics of the numbers
${\mathcal{N}}_t^{\times}$ of spontaneous returns of the walker to the origin
and ${\mathcal{N}}_t^{\bullet}$ of resetting events up to some observation time $t$.
These numbers are extensive in time in a strong sense:
all their joint cumulants grow linearly in $t$, with explicitly computable amplitudes,
and their fluctuations are described by a smooth bivariate large deviation function.
A non-trivial crossover phenomenon takes place in the regime of weak resetting and late times.
Remarkably, the time intervals between spontaneous returns to the origin of the reset random walk
form a renewal process described in terms of a single `dressed' probability distribution.
These time intervals are probabilistic copies of the first one, the `dressed' first-passage time.
The present work follows a broader study, covered in a companion paper,
on general nested renewal processes.
}

\maketitle

\section{Introduction}
\label{sec:intro}

This work builds upon a previous study on the replication of a renewal process
at random times, which is equivalent to nesting two generic renewal processes,
or, alternatively, to
considering a renewal process subject to random resetting~\cite{nested}.
In that study, we investigated the interplay between the two probability laws
governing the distribution of time intervals between renewals, on the one hand,
and resettings, on the other hand, resulting in a phase diagram that highlights
a rich range of behaviours.

\begin{figure}[!ht]
\begin{center}
\includegraphics[angle=0,width=0.7\linewidth,clip=true]{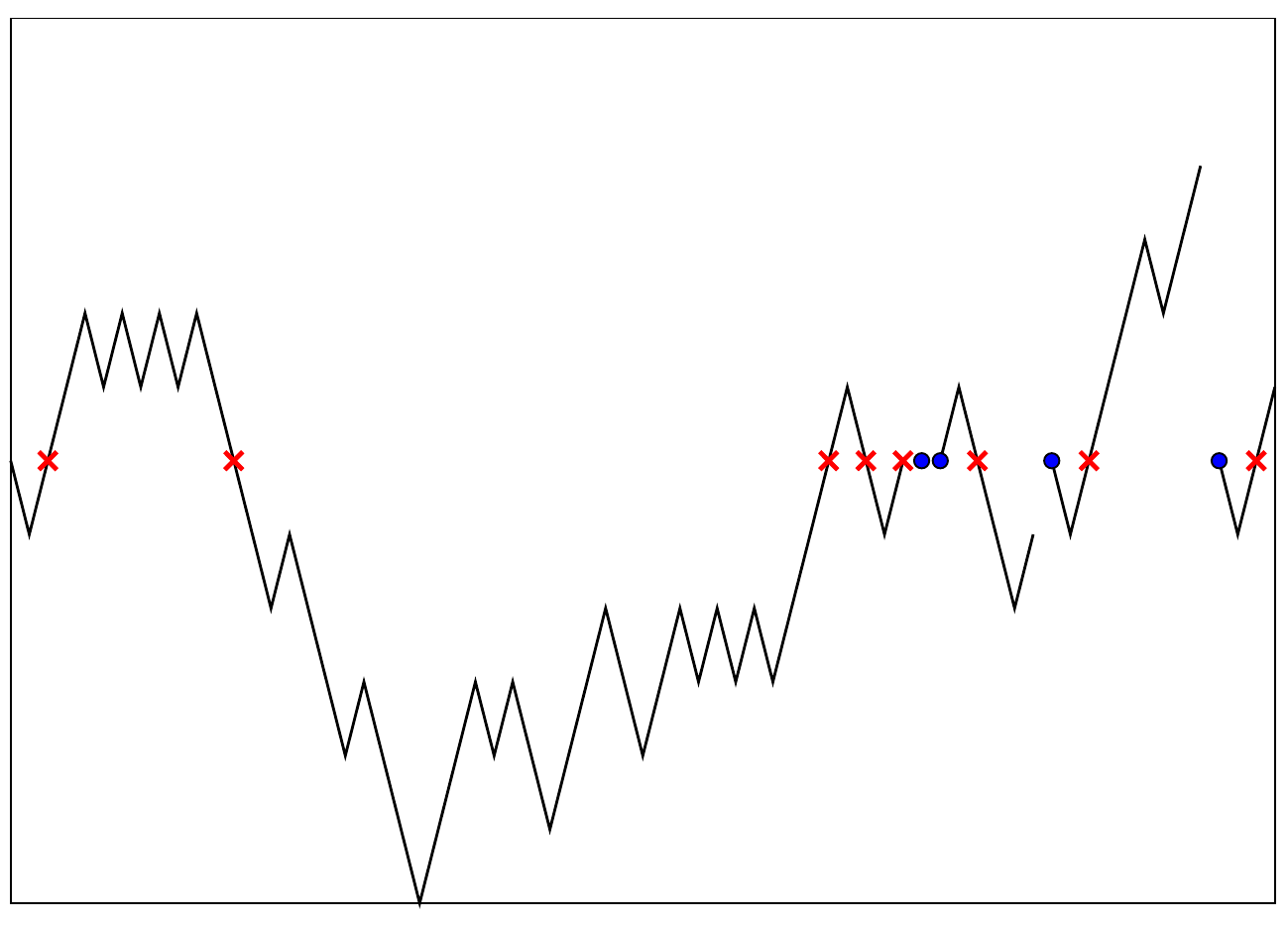}
\caption{\small
Example of a path of the P\'olya walk on the one-dimensional lattice under
stochastic resetting, generated by a simulation with $r=0.08$.
The walk starts at the origin.
It restarts afresh at the origin at each resetting event, figured by a dot.
Spontaneous returns to the origin are figured by crosses.}
\label{fig:snap}
\end{center}
\end{figure}

\begin{figure}[!ht]
\begin{center}
\includegraphics[angle=0,width=0.8\linewidth,clip=true]{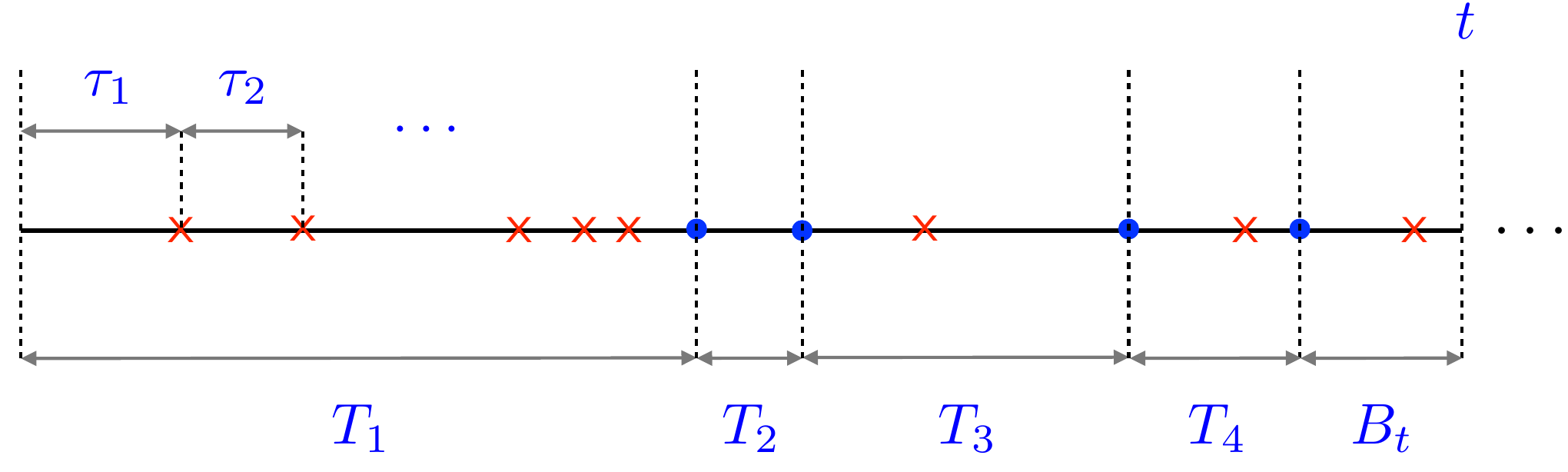}
\caption{\small
Sketch of the temporal events for the path of figure~\ref{fig:snap}.
Spontaneous returns to the origin of the walk are figured by crosses, resetting
events by dots.
The intervals of time between two crosses, $\ttau_1,\ttau_2,\dots$, have common
distribution $\rho(\tau)$ (see~(\ref{eq:defrho}), (\ref{eq:fctngen})).
The intervals of time between two resettings, $\T_1,\dots,\T_4$, have the
geometric distribution~(\ref{eq:geom}).
The last interval, $B_t$, represents the backward recurrence time, or age of
the
resetting process at time~$t$, i.e., the time elapsed since the previous
resetting event.}
\label{fig:nested}
\end{center}
\end{figure}

In the present work,
we investigate the specific case where the internal
renewal process consists of the epochs of returns to the origin
of the simple random walk (or P\'olya walk~\cite{polya}) on the one-dimensional
lattice,
while the external one involves discrete-time reset events
at which the process is restarted from the origin with probability~$r$ at each
time step.
The position $x_t$ of the walker at discrete time $t$ thus obeys the recursion
\beq\label{eq:defS}
x_{t+1}=\left\{
\begin{matrix}
0\hfill &\hbox{with probability $r$},\hfill\cr
x_t+\eta_{t+1}\quad &\hbox{with probability $1-r$},
\end{matrix}
\right.
\eeq
where $\eta_t=\pm1$ with equal probabilities.
The walk starts at the origin, $x_0=0$.
Figure~\ref{fig:snap} illustrates a sample path of the walk,
showing spontaneous returns to the origin marked by crosses and reset events
marked by dots.
Figure~\ref{fig:nested} provides a depiction of these temporal events and of
the intervals of time between them.

Much of the research in the theory of resetting processes has predominantly
concentrated on
continuous time stochastic processes (see~\cite{emsrev} for a review).
In contrast, relatively less emphasis has been devoted to discrete-time
processes.
An illustrative example of these processes involves discrete-time random walks
with continuous distributions of steps subject to resetting~\cite{kusmierz}.
Recent studies have delved into the statistics of extremes and records of
symmetric random walks with stochastic resetting~\cite{m2s2,glmax}.
Furthermore, investigations into discrete-time lattice random walks with
resetting have also been carried out.
Examples include unidirectional random walks with random
restarts~\cite{montero2016}, random walks where the walker is relocated to the
previous maximum~\cite{majumdar2015}, and random walks with preferential
relocations to previously visited locations~\cite{boyer}.

However, it is noteworthy that the P\'olya walk subject to resetting, defined
by (\ref{eq:defS}), has received relatively limited attention in the
literature~\cite{bonomo2021,bonomo2021-b,kumar,glmax,glsur}.
References~\cite{bonomo2021,bonomo2021-b} deal with general first-passage
properties of lattice random walks in discrete time, with application to the
P\'olya walk, while~\cite{kumar} contains a study of some aspects of the
statistics of records for the same walk.
In~\cite{glsur}, an analysis of the survival probability of symmetric random
walks with stochastic resetting was performed, specifically focussing on the
probability for the walker not to cross the origin up to time $t$, including
the example of the P\'olya walk~(\ref{eq:defS}).
Finally, the statistics of extremes and records for the P\'olya walk with
stochastic resetting are discussed in~\cite{glmax}.

The focus of the present work is on the joint statistics
of the numbers $\nc$ of spontaneous returns to the origin of the reset P\'olya
walk~(\ref{eq:defS}),
and $\nb$, denoting the count of reset events, up to a given time $t$.
These are the simplest observables one can think of for this process.
Their sum $\ncb=\nc+\nb$ is the total time spent by the walker at the origin.

The motivation for such a research stems from the analysis presented in the
companion paper~\cite{nested}.
The latter predicts that the more regular of two nested renewal processes
always governs the overall regularity of the entire process.
Here, the two renewal processes in question are made of the sequence of
spontaneous returns to the origin of the P\'olya walk, on the one hand, and of
the sequence of resetting events, on the other hand.
The latter---the more regular process of the two---is a Bernoulli process, as
can be seen on its definition~(\ref{eq:defS}).
In such a circumstance, as demonstrated in~\cite{nested}, $\mean{\nc}$ grows
linearly in time
and typical fluctuations of $\nc$ around its mean value are relatively
negligible.
The purpose of the present work is to corroborate these general results
and complete them by a thorough quantitative analysis
of the simple specific case at hand---the P\'olya walk under stochastic
resetting.

The setup and the main outcomes of this research are as follows.
Section~\ref{sec:basics} gives an exposition of background concepts and
results.
For the P\'olya walk without resetting (section~\ref{sec:pwwr}),
we recall results concerning the distribution of the intervals
between consecutive returns to the origin,
and the statistics of the number $N_t$ of such returns up to some time $t$.
Section~\ref{sre} contains a reminder on the statistics of resetting events
in discrete time.
Section~\ref{sec:obs} presents the detailed derivation
of the joint probability generating function of the random variables $\nc$ and
$\nb$ at any finite time $t$.
As a first application of the key equation~(\ref{eq:key}),
derived in section~\ref{sec:key},
the mean values of $\nc$ and $\nb$ are shown to grow linearly in time, as
\be
\mean{\nc}\approx\ac t,\qquad
\mean{\nb}\approx rt,\qquad
\mean{\ncb}\approx At,
\ee
where the amplitude
\be
A=\ac+r=\sqrt\frac{r}{2-r}
\ee
is identified with the steady-state probability for the walker to be at the origin.
In addition, we give in section~\ref{sec:dress}
an interpretation of the distribution of $\nc$ in terms of
a single `dressed' renewal process, and discuss its consequences.
An in-depth investigation of the statistics of $\nc$ and $\nb$ in the late-time regime
is done in section~\ref{sec:deviations},
highlighting the fact that these quantities are extensive in a strong sense.
First, their joint cumulants grow linearly in time, as
\be
\mean{(\nc)^k(\nb)^\ell}_c\approx c_{k,\ell}\,t.
\ee
We provide a method to evaluate all the cumulant amplitudes $c_{k,\ell}$,
and we give the explicit expressions of the first amplitudes
corresponding to $k+\ell\le3$ (see~(\ref{eq:c1})--(\ref{eq:c3})).
We also find that large fluctuations of $\nc$ and $\nb$ far from their mean values
obey a large deviation formula of the form
\be
\prob(\nc\approx\xi t,\ \nb\approx\eta t)\sim \e^{-I(\xi,\eta)t},
\ee
where the bivariate large deviation function $I(\xi,\eta)$
is the Legendre transform of the bivariate entropy function $S(\lam,\mu)$
generating the cumulant amplitudes $c_{k,\ell}$.
The ensuing univariate large deviation functions $\ib(\eta)$, $\ic(\xi)$, and
$I(\varphi)$, corresponding respectively to $\nb$, $\nc$ and their sum $\ncb$,
are plotted in figure~\ref{fig:lds}.
In the crossover regime at weak resetting and late times,
studied in section~\ref{sec:cross},
it is found that
\be
\nc\approx\sqrt{t}\,\mzeta,
\ee
where the rescaled random variable $\mzeta$ has a limiting distribution with
density $f(\zeta,u)$,
depending solely on the parameter $u=rt=\mean{\nb}$.
Figure~\ref{fig:fzeta} shows the density $f(\zeta,u)$
for several values of this parameter,
illustrating the crossover between a half-Gaussian form at $u=0$
and a drifting Gaussian at large $u$.
Section~\ref{sec:discuss} contains a brief discussion,
where the main outcomes of the present work
are put in perspective with those of the companion paper~\cite{nested}.
Some calculation details are relegated to two appendices.

\section{Background concepts}
\label{sec:basics}

\subsection{P\'olya walk without resetting}
\label{sec:pwwr}

As is well documented (see, e.g.,~\cite{feller1}), the sequence of returns to
the origin
of the P\'olya walk forms a discrete renewal process.
Let us denote by $\Tz$ the time of first return to the origin (from either
side),
and its distribution by
\beq\label{eq:defrho}
\rho(\tau)=\prob(\Tz=\tau).
\eeq
This quantity is non-zero whenever $\tau=2,4,\dots$ is an even integer.
This is also the common distribution of the intervals
between two consecutive returns to the origin,
denoted by $\ttau_1,\ttau_2,\dots$,
which are independent copies of $\Tz$.
The distribution $\rho(\tau)$ is known in terms of its generating function~\cite{feller1}
\beq\label{eq:fctngen}
\tilde\rho(z)=\sum_{\tau\ge0}z^\tau\rho(\tau)=1-\sqrt{1-z^2}.
\eeq
Introducing the binomial probabilities
\beq\label{eq:bdef}
b_n=\frac{(2n)!}{(2^nn!)^2}=\frac{{2n\choose n}}{2^{2n}},
\eeq
with generating function
\be
\tilde b(z)=\sum_{n\ge0}b_nz^n=\frac{1}{\sqrt{1-z}},
\ee
we have
\be
\rho(2n)=\frac{b_n}{2n-1}
\ee
for $n\ge1$, i.e.,
\beq\label{eq:rhofirst}
\rho(2)=\frac{1}{2},\quad\rho(4)=\frac{1}{8},
\quad\rho(6)=\frac{1}{16},\quad\rho(8)=\frac{5}{128},
\eeq
and so on.
When the even time $\tau$ becomes large, we have
\beq\label{eq:rhoasy}
\rho(\tau)
\approx\sqrt\frac{2}{\pi\tau^3}.
\eeq

The corresponding survival probability,
defined as the complementary distribution function of $\Tz$,
\beq\label{eq:defRtau}
R(\tau)=\prob(\Tz>\tau)=\sum_{j>\tau}\rho(j),
\eeq
obeys $R(\tau-1)-R(\tau)=\rho(\tau)$.
Its generating function reads
\beq\label{eq:rgen}
\tilde R(z)=\sum_{\tau\ge0}z^\tau R(\tau)=\frac{1-\tilde\rho(z)}{1-z}
=\frac{1+z}{\sqrt{1-z^2}}.
\eeq
We have therefore
\be
R(2n)=R(2n+1)=b_n,
\ee
i.e.,
\be
R(0)=R(1)=1,\quad R(2)=R(3)=\frac{1}{2},\quad R(4)=R(5)=\frac{3}{8},
\ee
and so on.
When $\tau$ becomes large, irrespective of its parity, we have
\beq\label{eq:rasy}
R(\tau)\approx\sqrt\frac{2}{\pi\tau}.
\eeq
The asymptotic estimate~(\ref{eq:rhoasy}) is minus twice the derivative
of~(\ref{eq:rasy}),
as it should be,
because~(\ref{eq:rhoasy}) only holds for even times $\tau$.

We now focus on the distribution of the number $N_t$
of returns of the walker to the origin up to time $t$.
This random variable is defined by the condition
\be
\ttau_1+\cdots+\ttau_{N_t}\le t<\ttau_1+\cdots+\ttau_{N_t+1},
\ee
hence the total time $t$ is decomposed into
\beq\label{eq:bt}
t=\ttau_1+\cdots+\ttau_{N_t}+b_t,
\eeq
where the last interval, $b_t$, is the backward recurrence time,
or the age of the renewal process at time $t$,
i.e., the elapsed time since the last return to the origin.
In the present discrete setting, $b_t=0,1,\dots,\ttau_{N_t+1}-1$.

A realisation of the set of random variables $\ttau_1,\dots,\ttau_{N_t},b_t$,
with $N_t=n$, denoted~by
\beq\label{eq:config}
\tilde\cc=\{\tau_1,\dots,\tau_{n},b\},
\eeq
has weight
\beq\label{eq:poids}
P(\tilde\cc)
=\rho(\tau_1)\ldots\rho(\tau_n)\,R(b)\,\delta\Big(\sum_{i=1}^n
\tau_i+b,t\Big),
\eeq
where $\delta(i,j)$ is the Kronecker delta symbol.

The distribution of $N_t$ ensues by summing the above weight
over all variables~$\{\tau_i\}$ and~$b$:
\be
p_n(t)=\prob(N_t=n)=\sum_{\{\tau_i\},b}\rho(\tau_1)\dots\rho(\tau_n) R(b)
\,\delta\Big(\sum_{i=1}^n\tau_i+b,t\Big).
\ee
The expression thus obtained is a discrete convolution, which is easier to handle
by taking its generating function with respect to $t$, which reads
\be
\sum_{t\ge0}w^t\,p_n(t)=\tilde\rho(w)^n\tilde R(w),
\ee
where $\tilde\rho(w)$ and $\tilde R(w)$
are respectively given by~(\ref{eq:fctngen}) and~(\ref{eq:rgen}).
The distribution of $N_t$ can be expressed compactly
through the probability generating function
\beq\label{eq:Zzt}
Z(z,t)=\mean{z^{N_t}}=\sum_{n\ge0}z^n\,p_n(t).
\eeq
The generating function of the latter quantity with respect to $t$ is
\be
\tilde Z(z,w)
=\sum_{t\ge0}w^tZ(z,t)
=\tilde R(w)\sum_{n\ge0}(z\tilde\rho(w))^n,
\ee
i.e.,
\beq\label{eq:Zzw}
\tilde Z(z,w)=\frac{1-\tilde\rho(w)}{(1-w)(1-z\tilde\rho(w))}.
\eeq

In particular,
the generating function with respect to $t$ of the mean number $\mean{N_t}$ of
returns reads
\beq\label{eq:fctngenN}
\sum_{t\ge0}w^t\mean{N_t}
=\frac{\partial}{\partial z}\tilde Z(z,w)\Big\vert_{z=1}
=\frac{\tilde\rho(w)}{(1-w)(1-\tilde\rho(w))}
=\frac{1+w}{(1-w^2)^{3/2}}-\frac{1}{1-w}.
\eeq
We have therefore
\be
\mean{N_{2n}}=\mean{N_{2n+1}}=(2n+1)b_n-1,
\ee
i.e.,
\be
\mean{N_0}=\mean{N_1}=0,\quad
\mean{N_2}=\mean{N_3}=\frac{1}{2},\quad
\mean{N_4}=\mean{N_5}=\frac{7}{8},
\ee
and so on.
When time $t$ becomes large, regardless of its parity, we have
\beq\label{eq:Ntasym}
\mean{N_t}\approx\sqrt\frac{2t}{\pi}.
\eeq
The probability of having $N_t=0$ is given by the generating function
\be
\sum_{t\ge0}w^t\,p_0(t)
=\tilde Z(0,w)
=\frac{1-\tilde\rho(w)}{1-w}
=\tilde R(w)
\ee
(see~(\ref{eq:rgen})).
We thus recover the expected result
\beq\label{eq:psur}
p_0(t)=\prob(\ttau>t)=R(t).
\eeq

The asymptotic distribution of $N_t$ in the regime of late times
can be extracted through a scaling analysis of~(\ref{eq:Zzw}).
Setting $w=\e^{-s}$ and $z=\e^{-p}$,
and working to leading order in the continuum regime where $s$ and $p$ are
small,
we obtain
\beq\label{eq:laplapn}
\int_0^\infty\dd t\,\e^{-st}\mean{\e^{-pN_t}}\approx\frac{1}{s+p\sqrt{s/2}}.
\eeq
Inverting the Laplace transforms in $p$ and in $s$ yields
\beq\label{eq:lapn}
\int_0^\infty\dd t\,\e^{-st}\,p_n(t)
\approx\sqrt\frac{2}{s}\,\e^{-\sqrt{2s}\,n},
\eeq
and finally
\beq\label{eq:hgau}
p_n(t)\approx\sqrt\frac{2}{\pi t}\,\e^{-n^2/(2t)}.
\eeq
We have thus recovered the known property that the asymptotic distribution
of the number~$N_t$ of returns to the origin of the simple random walk is a
half-Gaussian \cite{feller2}.
The limit of this distribution as $n\to0$
is consistent with the asymptotic behaviour of~$R(t)$ given by~(\ref{eq:rasy}).
The moments of the distribution~(\ref{eq:hgau}) read
\beq\label{eq:hgaumoms}
\mean{N_t^{2k}}\approx\frac{(2k)!}{2^kk!}\,t^k,\qquad
\mean{N_t^{2k+1}}\approx\sqrt\frac{2}{\pi}\,2^kk!\,t^{k+1/2}.
\eeq
In particular, the first moment agrees with~(\ref{eq:Ntasym}).

\subsection{Statistics of resetting events}
\label{sre}

The resetting events also constitute a discrete renewal process, referred to
in~\cite{nested} as the external renewal process.
The integer intervals of time $\T_1,\T_2\dots$ between two consecutive
resettings
have the geometric distribution
\beq\label{eq:geom}
f(T)=r(1-r)^{T-1}\qquad(T\ge1),
\eeq
whose complementary distribution function is given by
\be
\Phi(T)=\sum_{j>T}f(j)=(1-r)^T\qquad(T\ge0).
\ee
The corresponding generating functions read
\beqa
\tilde f(z)&=&\sum_{T\ge1}z^Tf(T)=\frac{rz}{1-(1-r)z},
\label{eq:serf}
\\
\tilde\Phi(z)&=&\sum_{T\ge0}z^T\Phi(T)=\frac{1-\tilde
f(z)}{1-z}=\frac{1}{1-(1-r)z}.
\label{eq:serFF}
\eeqa

The number of resetting events $M_t$ is defined by the condition
\be
\T_1+\cdots+\T_{M_t}\le t<\T_1+\cdots+\T_{M_t+1},
\ee
hence
\be
t=\T_1+\cdots+\T_{M_t}+B_t,\qquad B_t=0,1,\dots,\T_{M_t+1}-1.
\ee
The last interval $B_t$ is the backward recurrence time, or the age of the
resetting process at time $t$, i.e., the elapsed time since the last resetting
event.

A realisation of the set of random variables $\T_1,\dots,\T_{M_t},B_t$,
with $M_t=m$, denoted~by
\beq\label{eq:configR}
\cc=\{T_1,\dots,T_{m},B\},
\eeq
has weight
\beq\label{eq:poidsR}
P(\cc)
=f(T_1)\ldots f(T_m)\,\Phi(B)\;\delta\Big(\sum_{i=1}^m T_i+B,t\Big).
\eeq
Following the same approach as in~(\ref{eq:Zzw}), we have
\beq\label{eq:boules}
Y(y,w)
=\sum_{t\ge0}w^t\mean{y^{M_t}}
=\frac{\tilde\Phi(w)}{1-y\tilde f(w)}
=\frac{1}{1-(1-r+ry)w}.
\eeq
Hence
\beq\label{eq:ybin}
\mean{y^{M_t}}=(1-r+ry)^t.
\eeq
We thus recover that $M_t=0,\dots,t$ has the binomial distribution
\beq\label{eq:nbbin}
\prob(M_t=m)={t\choose m}r^{m}(1-r)^{t-m}
\eeq
at all times.
This known result agrees with the property that resetting events are independent from each other,
and therefore form a Bernoulli process.
In particular, the mean number of resettings reads
\beq\label{eq:meanMt}
\mean{M_t}=rt.
\eeq

\section{Spontaneous returns to the origin and resetting events}
\label{sec:obs}

\subsection{The key equation}
\label{sec:key}

As stated in the introduction, the main purpose of this work
is to analyse the joint distribution of the numbers $\nb$ of dots representing
resetting events and $\nc$ of crosses representing spontaneous returns to the
origin,
for the reset P\'olya walk up to time $t$.
These numbers are respectively given by
\be
\nb=M_t,
\label{nbmt}
\ee
introduced above, and
\beq\label{eq:Ntot+}
\nc=N_{\T_1-1}+N_{\T_2-1}+\cdots+N_{\T_{M_t-1}}+N_{B_t},
\eeq
as illustrated on figure~\ref{fig:snap}.
The sum of these two numbers is denoted by
\be
\ncb=\nc+\nb,
\ee
and reads
\beq\label{eq:Ntot++}
\ncb
=\sum_{\tau=1}^t\delta(x_\tau,0),
\eeq
where $x_\tau$ is the position of the walker at time $\tau$
(see~(\ref{eq:defS})).
This is the total time spent by the walker the origin,
either by a resetting event (a dot) or by a spontaneous return (a cross).

Notice that the distribution of the number of dots, $\nb=M_t$
(see~(\ref{nbmt})),
is known from~(\ref{eq:boules}), (\ref{eq:ybin}), (\ref{eq:nbbin})
and~(\ref{eq:meanMt}).
Thus
\beq\label{eq:boules+}
\Y(y,w)=\sum_{t\ge0}w^t\mean{y^{\nb}}=Y(y,w)=\frac{1}{1-(1-r+ry)w},
\eeq
\beq\label{eq:ybin+}
\mean{y^{\nb}}=(1-r+ry)^t,
\eeq
\beq\label{eq:nbbin+}
\prob(\nb=\nn)={t\choose\nn}r^\nn(1-r)^{t-\nn},
\eeq
\beq\label{eq:meannb}
\mean{\nb}=rt.
\eeq
Notice also that the situation simplifies in the following two special cases.
First,
in the absence of resetting ($r=0$), we have
\be
\nc=\ncb=N_t,\qquad\nb=0,
\ee
second, when a resetting event occurs at every time step ($r=1$), then
\beq\label{eq:run}
\nb=\ncb=t,\qquad\nc=0.
\eeq

The central quantity for the determination of the joint statistics of $\nc$ and
$\nb$ is
the generating function
\be
\Z(z,y,t)=\mean{z^{\nc}y^{\nb}},
\ee
where the average is taken over the external
configurations $\cc$ (see~(\ref{eq:configR})), with weight $P(\cc)$ given by
(\ref{eq:poidsR}), and over the internal configurations $\tilde\cc$
(see~(\ref{eq:config})), with weight
$P(\tilde\cc)$ given by (\ref{eq:poids}).
Thus
\be
\Z(z,y,t)=\sum_{\cc}P(\cc)\sum_{\tilde\cc}P(\tilde\cc)\,z^{\nc}y^{\nb},
\ee
with the notations
\be
\sum_{\cc}=\sum_{m\ge0}\sum_{\{T_i,\,B\}},
\qquad \sum_{\tilde\cc}=\sum_{n\ge0}\sum _{\{\tau_i,\,b\}}.
\ee
The average over the internal variables of each term $z^{N_{T_i}-1}$ with
weight $P(\tilde\cc)$ gives a factor $Z(z,T_i-1)$ (see~(\ref{eq:Zzt})).
We then average over the external variables with weight $P(\cc)$ to arrive at
\beq\label{eq:calN}
\mean{z^{\nc}y^{\nb}}=\sum_{\cc}P(\cc)y^mZ(z,T_1-1)\dots Z(z,T_m-1)Z(z,B).
\eeq
The expression thus obtained is a discrete convolution, which is easier to
handle
by taking its generating function with respect to $t$, leading to
\beq\label{eq:Zcal}
\tilde\Z(z,y,w)
=\sum_{t\ge0}w^t\Z(z,y,t)
=\sum_{m\ge0}y^m\tilde\varphi(z,w)^m\tilde\psi(z,w)
=\frac{\tilde\psi(z,w)}{1-y\tilde\varphi(z,w)},
\eeq
with
\bea
\tilde\varphi(z,w)
&=&\sum_{T\ge1}w^Tf(T)\,Z(z,T-1),
\label{phiex}
\\
\tilde\psi(z,w)
&=&\sum_{B\ge0}w^B\Phi(B)\,Z(z,B).
\label{psiex}
\eea

The expressions~(\ref{eq:Zcal})--(\ref{psiex}) are quite general~\cite{nested}.
They hold for arbitrary distributions $\rho(\tau)$ and $f(T)$,
both in a continuous and in a discrete setting.
For the case at hand,
the generating functions of $Z(z,T)$, $f(T)$ and $\Phi(T)$
are respectively given in~(\ref{eq:Zzw}),~(\ref{eq:serf}) and~(\ref{eq:serFF}).
We have therefore
\bea
\tilde\varphi(z,w)
&=&rw\,\tilde Z(z,\rw),
\nonumber
\\
\tilde\psi(z,w)
&=&\tilde Z(z,\rw),
\eea
where we introduced the shorthand notation
\beq\label{eq:rwdef}
\rw=(1-r)w,
\eeq
and finally
\beq\label{eq:key}
\tilde\Z(z,y,w)=\frac{\tilde Z(z,\rw)}{1-r yw\tilde Z(z,\rw)}.
\eeq
The expression~(\ref{eq:key}),
where $\tilde Z(z,w)$ is given in~(\ref{eq:Zzw}),
and $\tilde\rho(w)$ in~(\ref{eq:fctngen}),
is the key result of this section.

The distribution of $\nb=M_t$ is obtained by setting $z=1$ in~(\ref{eq:key}).
We thus recover the expression~(\ref{eq:boules+}),
since $\tilde\Z(1,y,w)=\Y(y,w)$, as it should be.
The distribution of~$\nc$ is obtained by setting $y=1$ in~(\ref{eq:key}).
We thus have
\beq\label{eq:Zz1w}
\tilde\Z(z,1,w)=\sum_{t\ge0}w^t\mean{z^{\nc}}
=\frac{\tilde Z(z,\rw)}{1-rw\tilde Z(z,\rw)}.
\eeq
This formula relates the same quantity with and without resetting.
Rational relationships of this form are specific to Poissonian resetting in
continuous time
or to geometric resetting in discrete time
(see, e.g.,~\cite{emsrev,kusmierz,m2s2,glmax,glsur,bonomo2021-b}).
Another interpretation of~(\ref{eq:Zz1w}) will be given in
section~\ref{sec:dress}.

The general expression of the mean value $\mean{\nc}$
ensues from~(\ref{eq:Zcal}), and reads
\bea
\sum_{t\ge0}w^t\mean{\nc}
&=&\frac{\dd}{\dd z}\Z(z,y,w)\Big\vert_{z=y=1}
\\
&=&\frac{1}{1-\tilde f(w)}\Bigg(\sum_{T\ge1}
w^{T}\mean{N_{T-1}}\frac{f(T)}{1-w}+\sum_{B\ge0} w^{B}\mean{N_B}\Phi(B)\Bigg).
\eea
In the case at hand, this expression yields
\beq\label{eq:FGNx}
\sum_{t\ge0} w^t\mean{\nc}
=\frac{(1-\rw)\tilde\rho(\rw)}{(1-w)^2(1-\tilde\rho(\rw))},
\eeq
which can alternatively be obtained from~(\ref{eq:key}).

In the absence of resetting ($r=0$),~(\ref{eq:FGNx}) gives
back~(\ref{eq:fctngenN}),
as it should be, since $\nc=N_t$.
In the presence of resetting ($r\ne0$),
we obtain the linear growth law
\beq\label{eq:ncasy}
\mean{\nc}\approx\ac t
\eeq
in the late-time regime,
with
\beq\label{eq:acres}
\ac=\frac{r\tilde\rho(1-r)}{1-\tilde\rho(1-r)}=\sqrt\frac{r}{2-r}-r.
\eeq
In the regime of weak resetting,
the mean value $\mean{\nc}$ exhibits
a smooth crossover between the square-root law~(\ref{eq:Ntasym})
and the linear law~(\ref{eq:ncasy}).
The complete determination of the distribution of~$\nc$ throughout this
crossover regime
will be given in section~\ref{sec:cross}.

\begin{figure}
\begin{center}
\includegraphics[angle=0,width=0.7\linewidth,clip=true]{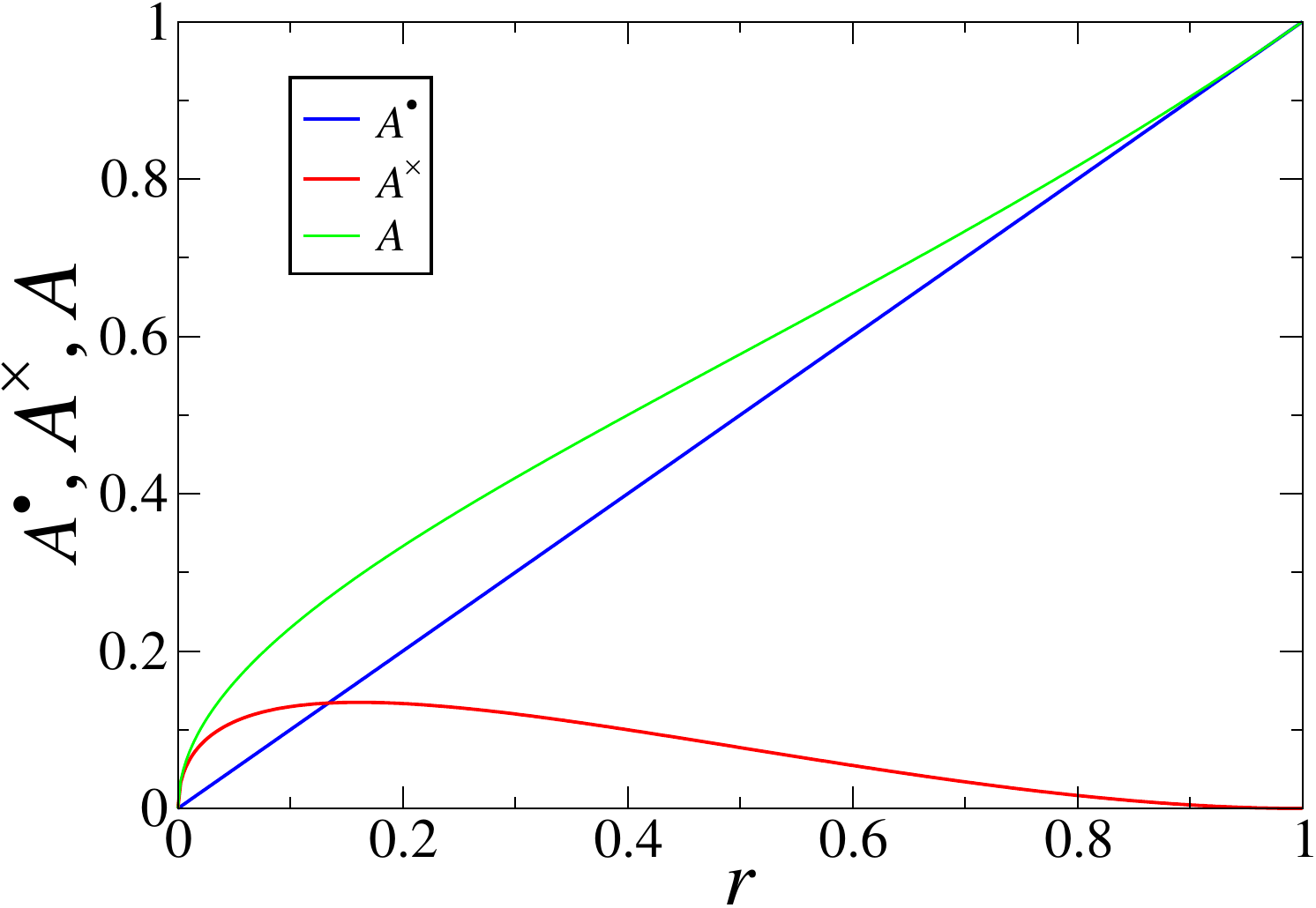}
\caption{\small
Amplitudes $\ab=r$, $\ac$ and their sum $A$
entering the growth laws~(\ref{eq:meannb}),~(\ref{eq:ncasy}) and~(\ref{ncbasy})
of $\mean{\nb}$, $\mean{\nc}$ and $\mean{\ncb}$,
plotted against the resetting probability $r$.}
\label{fig:amps}
\end{center}
\end{figure}

Summing expressions~(\ref{eq:meannb}) and~(\ref{eq:ncasy}),
we end up with
\be
\mean{\ncb}\approx At,
\label{ncbasy}
\ee
with
\beq\label{eq:ares}
A=\sqrt\frac{r}{2-r}.
\eeq
Both amplitudes $\ac$ and $A$ vanish as
\beq\label{eq:Ax0}
\ac\approx A\approx\sqrt\frac{r}{2}
\eeq
as $r\to0$.
This square-root scaling will be corroborated by the analysis of the crossover
regime (see section~\ref{sec:cross}).
The amplitude $\ac$ vanishes quadratically as
\beq\label{eq:Ax1}
\ac\approx\frac{(1-r)^2}{2}
\eeq
as $r\to1$,
testifying that the presence of a cross,
i.e., of a spontaneous return of the walker to the origin,
requires at least two successive time steps without resetting.
This amplitude reaches its maximum $\ac_{\max}=0.134884\dots$ for
$r=0.160713\dots$
As for the amplitude $A$ of $\mean{\ncb}$, it increases monotonically from 0 to
1
as $r$ increases in the same range of values.
Figure~\ref{fig:amps} shows plots of the amplitudes $\ab=r$, $\ac$,
and of their sum $A$, against the resetting probability $r$.

The formula~(\ref{eq:ares}) can be compared with the expression
for the position distribution of the walker in the nonequilibrium stationary state
reached in the limit of infinitely large times, namely
\be
p(x)=\sqrt\frac{r}{2-r}\,\lam^{-\abs{x}},\qquad\lam=\frac{1+\sqrt{r(2-r)}}{1-r}.
\label{pst}
\ee
For the sake of completeness,
we give a brief self-contained derivation of this result in Appendix~\ref{sec:appb}.
It had been derived in a general setting in~\cite{glmax}.
The distribution~(\ref{pst}) falls off exponentially with the distance
to the resetting point, i.e., the origin,
where it reaches its maximum
\be
p(0)=\sqrt\frac{r}{2-r}.
\ee
The identity $A=p(0)$ is to be expected,
as both sides represent the fraction of time spent by the walker at the origin.
Formally, this identity can be derived
by taking the mean values of~(\ref{eq:Ntot++})
and using the fact that the occupation probability of the origin at time $t$,
$\mean{\delta(x_t,0)}$, tends to $p(0)$ at late times.

\subsection{The reset process seen as a `dressed' renewal process}
\label{sec:dress}

\begin{figure}
\begin{center}
\includegraphics[angle=0,width=0.8\linewidth,clip=true]{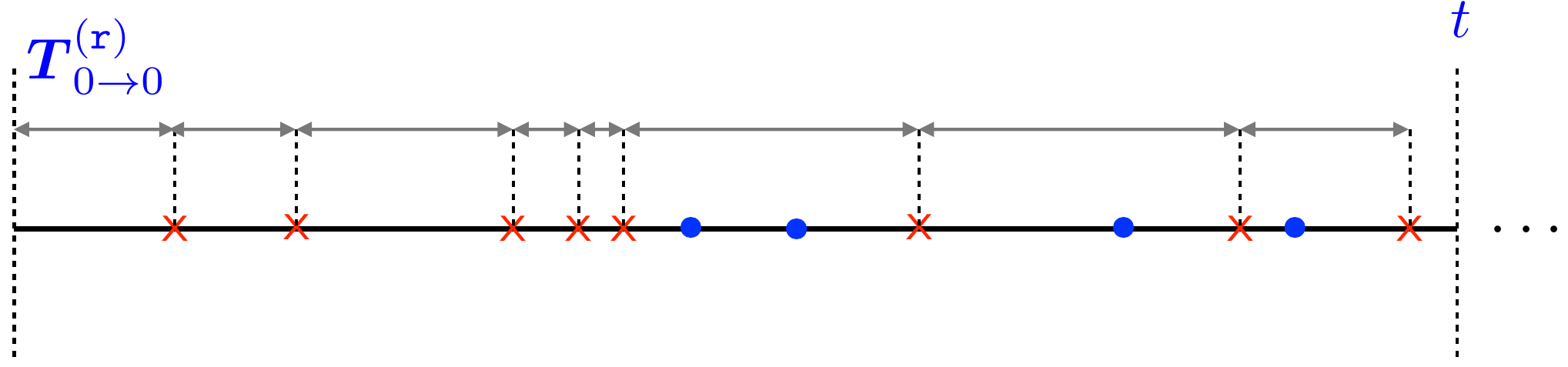}
\caption{\small
The time intervals between spontaneous returns to the origin (crosses)
of the reset random walk form a renewal process,
described in terms of the `dressed' probability
distribution~(\ref{eq:dressed2}).
These time intervals are probabilistic copies of the first one,
$\Tz^{(\tt r)}$, the `dressed' first-passage time.
(Compare to figure~\ref{fig:nested}.)}
\label{fig:dressed}
\end{center}
\end{figure}

It is interesting to note that there exists another interpretation
of the expression~(\ref{eq:Zz1w}), encoding the full statistics of $\nc$.
This expression is of the form~(\ref{eq:Zzw}),
up to the replacement of $\tilde\rho(w)$ by
\beq\label{eq:dressed}
\tilde\rho^{(\tt r)}(w)=\frac{(1-\rw)\tilde\rho(\rw)}{1-w+rw\tilde\rho(\rw)}.
\eeq
This implies that the interarrival times between spontaneous returns to the origin
(crosses in figure~\ref{fig:nested})
remain, in the presence of resetting, independent, identically distributed,
random variables, whose common probability distribution,
\beq\label{eq:dressed2}
\rho^{(\tt r)}(\tau)=\prob(\Tz^{(\tt r)}=\tau),
\eeq
has a generating function given by~(\ref{eq:dressed}),
and depends only on the resetting probability~$r$ (see~(\ref{eq:firstval})).
In other words, these interarrival times form a renewal process,
defined by the `dressed' distribution~(\ref{eq:dressed2}),
i.e., they are probabilistic replicas of the `dressed' first-passage time $\Tz^{(\tt r)}$
(or, in the present context, time of first return to the origin),
as depicted in figure~\ref{fig:dressed}.
The superscript in these expressions is an abbreviation
for \textit{replication} or \textit{resetting}.

The first few values of $\rho^{(\tt r)}(\tau)$ read
\beqa\label{eq:firstval}
\rho^{(\tt r)}(1)&=&0,\quad \rho^{(\tt r)}(2)=\frac{1}{2}(1-r)^2,\quad
\rho^{(\tt r)}(3)=\frac{1}{2}r(1-r)^2,\quad
\nonumber\\
\rho^{(\tt
r)}(4)&=&\frac{1}{8}(1-r)^4+\frac{1}{2}r^2(1-r)^2+\frac{1}{2}r(1-r)^3=\frac{1}{8}(1-r^2)^2.
\eeqa
The first three expressions are easy to guess,
whereas the fourth is the sum of the probabilities of the following events:
$\{$return to the origin in four steps without any resetting$\}$, $\{$two
resettings first, then return to the origin in two steps$\}$, and finally
$\{$one step away from the origin, a resetting, then return to the origin in two steps$\}$.
The right-hand sides in~(\ref{eq:firstval}) give back~(\ref{eq:rhofirst}), when $r=0$.

The existence of the aforementioned dressed renewal process is a remarkable phenomenon,
which occurs whether the time intervals between resetting events,
pertaining to the external renewal process,
follow an exponential distribution in continuous time (resulting in a Poisson process),
as described in~\cite{nested},
or a geometric distribution (see~(\ref{eq:geom})) in discrete time
(resulting in a Bernoulli process), as described above.
This allows us to access several characteristic features of the reset P\'olya walk
by considering quantities pertaining to the dressed renewal process, as we now elaborate.

The `dressed' survival probability is naturally defined, in line with~(\ref{eq:defRtau}), as
\beq\label{eq:dressedR}
R^{(\tt r)}(\tau)=\prob(\Tz^{(\tt r)}>\tau)=\sum_{j>\tau}\rho^{(\tt r)}(j).
\eeq
Starting from~(\ref{eq:dressed}) and using~(\ref{eq:rgen}),
as well as the corresponding generating function for the dressed distribution,
\be
\tilde R^{(\tt r)}(w)=\frac{1-\tilde \rho^{(\tt r)}(w)}{1-w},
\ee
we can establish the following connection between the generating functions
for the survival probabilities~(\ref{eq:defRtau}) and~(\ref{eq:dressedR})
in the absence or in the presence of resetting,
\beq\label{eq:kus}
\tilde R^{(\tt r)}(w)=\frac{\tilde R(\rw)}{1-rw \tilde R(\rw)}.
\eeq
The same relation holds for the probability of not crossing the origin up to
integer time $t$,
for a discrete-time walker with continuous steps~\cite{kusmierz},
or for the P\'olya walk~\cite{glsur}.

The moments of $\Tz^{(\tt r)}$ can be derived from~(\ref{eq:dressed}).
We have in particular
\be
\mean{\Tz^{(\tt r)}}=\frac{1}{\ac},
\label{acave}
\ee
meaning that $\mean{\Tz^{(\tt r)}}$ has a minimum at $r=0.160713\dots$,
and that~(\ref{eq:ncasy}) can be recast as
\be
\mean{\nc}\approx\frac{t}{\mean{\Tz^{(\tt r)}}},
\ee
which is consistent with intuition.

Another quantity of interest is the probability of occurrence of a cross
(spontaneous return to the origin of the P\'olya walk) at time $t$, in the absence or in the presence of resetting, i.e., respectively,
\beq\label{eq:udef}
U(t)=\mean{N_t}-\mean{N_{t-1}},\qquad
U^{(\tt r)}(t)=\mean{\nc}-\mean{\ncminus}\qquad (t\ge1),
\eeq
completed by $U(0)=U^{(\tt r)}(0)=1$.
The corresponding generating functions are given by (see~(\ref{eq:fctngenN}))
\be
\tilde U(w)=\frac{1}{1-\tilde\rho(w)},\qquad
\tilde U^{(\tt r)}(w)=\frac{1}{1-\tilde\rho^{(\tt r)}(w)}.
\ee
We have $U(2n)=b_n$ (see~(\ref{eq:bdef})), whereas $U(2n+1)=0$.
For $t$ large, $U^{(\tt r)}(t)$ converges very rapidly to $\ac$.
The expression of $U^{(\tt r)}(t)$ in the crossover regime of weak resetting and late times will be given in~(\ref{usca}).

\begin{figure}
\begin{center}
\includegraphics[angle=0,width=0.7\linewidth,clip=true]{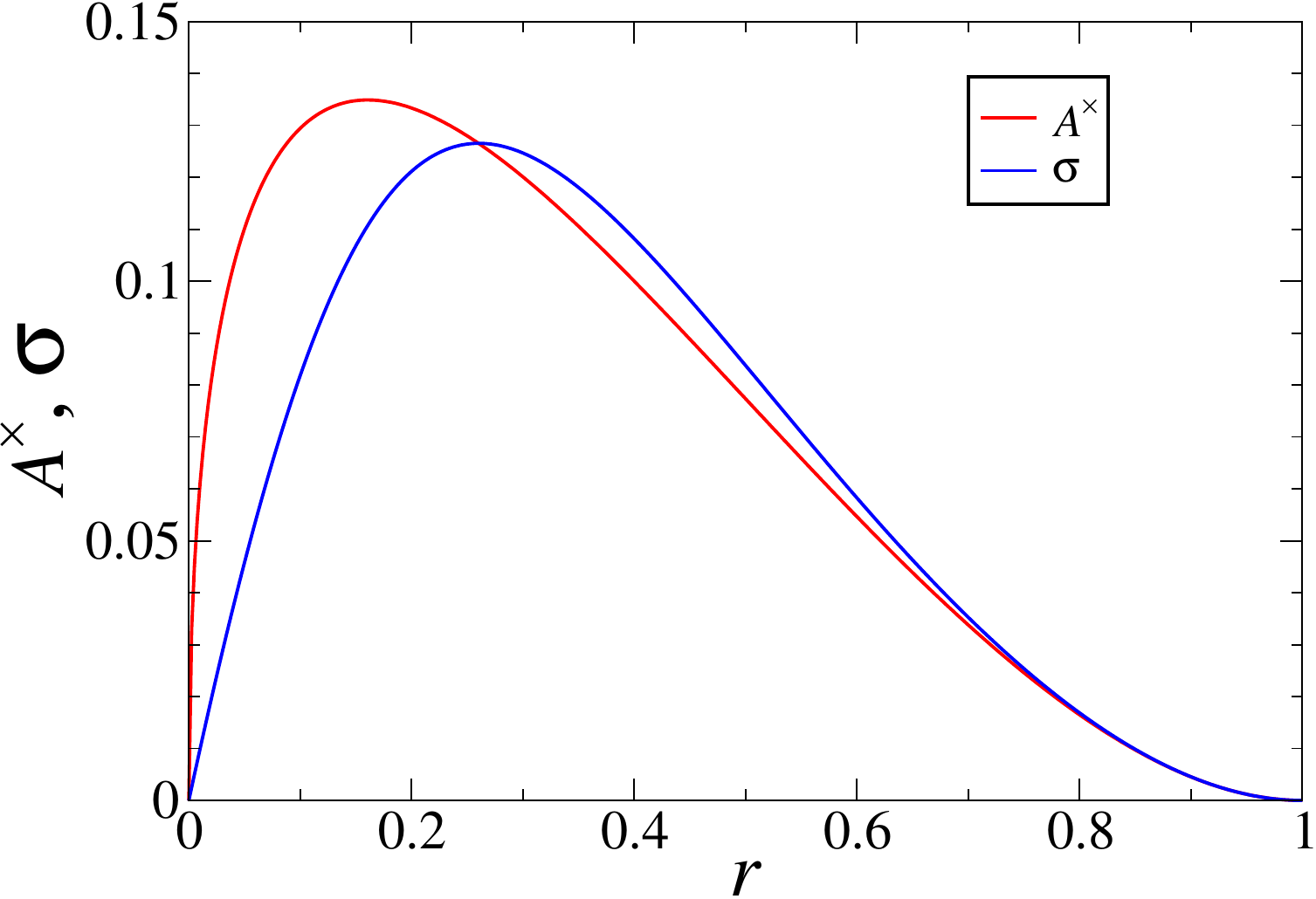}
\caption{\small
Amplitude $\ac$ (see~(\ref{eq:ncasy}) and~(\ref{acave}))
and decay rate $\sigma$ (see~(\ref{eq:sigrat})),
plotted against the resetting probability $r$.}
\label{fig:asigma}
\end{center}
\end{figure}

The tail of the dressed distribution $\rho^{(\tt r)}(\tau)$
falls off exponentially as
\beq\label{eq:sigrat}
\rho^{(\tt r)}(\tau)\sim\e^{-\sigma\tau}.
\eeq
The decay rate $\sigma$ is such that $w_0=\e^\sigma$
is the smallest zero of the denominator of~(\ref{eq:dressed}), obeying
\be
r^2(1-r)w_0^3+r^2w_0^2+(1-r)w_0-1=0.
\ee
The dependence of the decay rate $\sigma$ on $r$
is qualitatively similar to that of the amplitude~$\ac$
(see~figure~\ref{fig:asigma}).
It vanishes as $\sigma\approx r$ as $r\to0$,
and as $\sigma\approx (1-r)^2/2$
as $r\to1$ (compare to~(\ref{eq:Ax0}) and~(\ref{eq:Ax1})),
reaching its maximum $\sigma_{\max}=0.126530\dots$ for $r=0.260465\dots$

As a summary, while the original renewal process made of the spontaneous
returns to the origin of the P\'olya walk is not stationary,
since the distribution~(\ref{eq:rhoasy}) has a fat power-law tail,
the dressed renewal process defined by the exponentially decaying
distribution~(\ref{eq:sigrat}) becomes eventually stationary (see~\cite{glrenew,nested}).

The formula~(\ref{eq:dressed}) is presented in~\cite{bonomo2021-b} within a different context.
It appears as an application to geometric resetting of a general formula
concerning the distribution of first-passage times in a generic discrete-step
walk in the presence of an arbitrary distribution of reset events.
This reference solely focuses on first-passage properties.
In particular, it does not delve into the renewal structure of the spontaneous
returns to the origin of the reset P\'olya walk defined by~(\ref{eq:defS}).
The method used in~\cite{bonomo2021-b} to derive the distribution of the first-passage time
for an arbitrary distribution of reset events, which builds upon prior
research~\cite{reuveniPRL,reuveni,bonomo2021}, is revisited in~\cite{nested}.

\section{Cumulants and large deviations in the late-time regime}
\label{sec:deviations}

In this section we continue the investigation of the joint statistics of the
quantities $\nc$
(number of dots) and $\nb$ (number of crosses) for the P\'olya walk
in the late-time regime.
We shall demonstrate that these variables are extensive in a strong sense,
first by examining their joint cumulants
and then by investigating the corresponding large deviation functions.

\subsection{Cumulants}

The starting point of the analysis is again the key formula~(\ref{eq:key}).
The late-time regime is governed by the smallest zero $w_\star(z,y)$
of the denominator of that formula, which entails an exponential law of the
form
\beq\label{eq:exes}
\mean{z^{\nc}y^{\nb}}\sim w_\star(z,y)^{-t}
\eeq
for the joint probability generating function of $\nc$ and $\nb$ in the
late-time regime.
Introducing the notations
\beq\label{eq:lmnots}
z=\e^\lam,\qquad y=\e^\mu,
\eeq
and
\beq\label{eq:lmnots+}
w_\star(z,y)=\e^{-S(\lam,\mu)}
\eeq
brings the estimate~(\ref{eq:exes})
for the generating function of the joint cumulants of $\nc$ and $\nb$
to the more familiar form
\beq\label{eq:s2}
\mean{\e^{\lam\nc+\mu\nb}}\sim\e^{S(\lam,\mu)t}.
\eeq
The exponential law~(\ref{eq:s2})
implies that all joint cumulants grow linearly with time, as
\be
\mean{(\nc)^k(\nb)^\ell}_c\approx c_{k,\ell}\,t,
\ee
where the amplitudes $c_{k,\ell}$ are the coefficients of the series expansion
\beq\label{eq:sers}
S(\lam,\mu)=\sum_{k+\ell\ge1}c_{k,\ell}\,\frac{\lam^k}{k!}\,\frac{\mu^\ell}{\ell!}
\eeq
of the entropy function $S(\lam,\mu)$ entering~(\ref{eq:s2}).

In particular, setting $\lam=\mu$ in~(\ref{eq:s2}) yields
\be
\mean{\e^{\lam\ncb}}\sim\e^{S(\lam,\lam)t},
\ee
so that $S(\lam,\lam)$
generates the amplitudes of the cumulants of $\ncb=\nc+\nb$ in the late-time
regime.
We have
\be
\mean{(\ncb)^n}_c\approx C_n\,t,
\ee
with
\be
S(\lam,\lam)=\sum_{n\ge1}C_n\,\frac{\lam^n}{n!}
\ee
and
\beq\label{eq:csum}
C_n=\sum_{k=0}^n{n\choose k}c_{k,n-k}.
\eeq

In order to derive explicit expressions for the function $S(\lambda,\mu)$
defined
in~(\ref{eq:lmnots+}), and therefore for the cumulant amplitudes $c_{k,\ell}$
and $C_n$,
we need an expression of the smallest zero of the denominator
in~(\ref{eq:key}).
This is done in Appendix~\ref{sec:appa},
yielding $w_\star(z,y)=w_1$, where $w_1$ is known explicitly
from~(\ref{eq:lmnots}),~(\ref{eq:pdefs}),~(\ref{eq:bpq}),
(\ref{eq:wres}) and~(\ref{eq:sigth}), whence, finally,
\beq\label{eq:s2trig}
S(\lam,\mu)=-\ln w_1.
\eeq

By expanding the expression above as a power series in $\lam$ and $\mu$,
we obtain explicit expressions for the cumulants $c_{k,\ell}$ and $C_n$,
which can be further reduced to expressions linear~in
\be
A=\sqrt\frac{r}{2-r}
\ee
(see~(\ref{eq:ares})),
with coefficients rational in $r$.
The first few formulas given below
testify that their complexity increases very fast with the order of the
cumulants.
To first order in $\lam$ and $\mu$, we have
\beqa
c_{1,0}&=&A-r,
\nonumber\\
c_{0,1}&=&r,
\nonumber\\
C_1&=&A,
\label{eq:c1}
\eeqa
in agreement with~(\ref{eq:meannb}),~(\ref{eq:acres}) and~(\ref{eq:ares}).
To second order, we have
\beqa
c_{2,0}&=&-\frac{4-3r}{2-r}\,A+\frac{2+4r-9r^2+5r^3-r^4}{(2-r)^2},
\nonumber\\
c_{1,1}&=&\frac{1-r}{2-r}\,A-r(1-r),
\nonumber\\
c_{0,2}&=&r(1-r),
\nonumber\\
C_2&=&-A+\frac{2-r^2}{(2-r)^2}.
\label{eq:c2}
\eeqa
To third order, we have
\beqa
c_{3,0}&=&\frac{3+38r-76r^2+46r^3-8r^4}{r(2-r)^3}\,A
\nonumber\\
&-&\frac{12+14r-66r^2+73r^3-43r^4+15r^5-2r^6}{(2-r)^3},
\nonumber\\
c_{2,1}&=&-\frac{(1-r)(4-5r)}{(2-r)^2}\,A
\nonumber\\
&+&\frac{r(1-r)(12-32r+30r^2-13r^3+2r^4)}{(2-r)^3},
\nonumber\\
c_{1,2}&=&\frac{(1-r)(1-2r)}{(2-r)^2}\,A-r(1-r)(1-2r),
\nonumber\\
c_{0,3}&=&r(1-r)(1-2r),
\nonumber\\
C_3&=&\frac{3+20r-31r^2+10r^3+r^4}{r(2-r)^3}\,A-\frac{3(2-r^2)}{(2-r)^2}.
\label{eq:c3}
\eeqa

More specific results can be derived in a few special and limiting situations.
We have observed that the statistics of $\nb$ is simple,
as its distribution is the binomial~(\ref{eq:nbbin+}).
This leads to two consequences.
First, all its cumulants have an exact linear behaviour in $t$ at all times:
\be
\mean{(\nb)^\ell}_c=c_{0,\ell}\,t.
\ee
Second, the cumulant amplitudes $c_{0,\ell}$ can be derived
from~(\ref{eq:ybin+}),
which amounts~to
\beq\label{eq:smu}
S(0,\mu)=\sum_{\ell\ge1}c_{0,\ell}\,\frac{\mu^\ell}{\ell!}=\ln(1-r+r\e^\mu).
\eeq
The first few amplitudes read
\bea
c_{0,1}&=&r,\qquad
c_{0,2}=r(1-r),
\nonumber\\
c_{0,3}&=&r(1-r)(1-2r),\qquad
c_{0,4}=r(1-r)(1-6r+6r^2).
\eea
The first three expressions agree with~(\ref{eq:c1})--(\ref{eq:c3}).
The amplitudes $c_{0,\ell}$ are polynomials in $r$ of increasing degrees,
obeying the linear differential recursion~\cite{haldane}

\be
c_{0,\ell+1}=r(1-r)\frac{\dd c_{0,\ell}}{\dd r}\qquad(\ell\ge1).
\ee
They read explicitly
\beq\label{eq:cstir}
c_{0,\ell}=\sum_{k=1}^\ell(-1)^{k-1}(k-1)!\stwo{\ell}{k}r^k,
\eeq
where $\stwo{\ell}{k}$ are the Stirling numbers of the second kind.

In the regime of strong resetting $(r\to1)$, we have
\be
S(\lam,\mu)
=\mu+(1-r)(\e^{-\mu}-1)
+(1-r)^2(\e^{-\mu}-\e^{-2\mu}+\half(\e^{\lam-2\mu}-1))+\cdots
\ee
For $r=1$, only $c_{0,1}=1$ is non-zero, in agreement with~(\ref{eq:run}).
To first order in $1-r$, only the $c_{0,\ell}$ are non-zero,
as there are no crosses at this order.
Their behaviour $c_{0,\ell}\approx(-1)^\ell(1-r)$ for $\ell\ge2$ agrees
with~(\ref{eq:cstir}).
All cumulant amplitudes $c_{k,\ell}$ become non-trivial to second order in
$1-r$.

The regime of weak resetting $(r\to0)$ is more subtle.
This richness is related to the crossover phenomenon
that will be examined in section~\ref{sec:cross}.
An inspection of the formulas~(\ref{eq:c1}),~(\ref{eq:c2}) and~(\ref{eq:c3})
yields the scaling $c_{k,\ell}\sim A^{2-k}$, with $A\approx\sqrt{r/2}$
(see~(\ref{eq:ares})),
at least for odd $k$.
This observation can be corroborated by the following scaling analysis.
Let us assume that the cumulant amplitudes behave as
\beq\label{eq:csca}
c_{k,\ell}\approx b_{k,\ell}\,A^{2-k}
\eeq
as $A\to0$, i.e., $r\to0$,
where the $b_{k,\ell}$ are constants to be determined.
This ansatz translates into the scaling form
\beq\label{eq:ssca}
S(\lam,\mu)\approx A^2\,F(h,\mu),
\eeq
with
\be
h=\frac{\lam}{A},\qquad
F(h,\mu)=\sum_{k+\ell\ge1}b_{k,\ell}\,\frac{h^k}{k!}\,\frac{\mu^\ell}{\ell!}.
\ee
Inserting the scaling form~(\ref{eq:ssca}), with the
notations~(\ref{eq:lmnots}), into~(\ref{eq:wpol}),
and expanding to the first non-trivial order in $A$,
we obtain that the scaling function $F(h,\mu)$ obeys the quadratic equation
\be
2F^2-(8(\e^\mu-1)+h^2)F+8(\e^\mu-1)^2-2h^2=0,
\ee
hence
\beq\label{eq:fres}
F(h,\mu)=2(\e^\mu-1)+\frac{h^2}{4}+h\sqrt{\e^\mu+\frac{h^2}{16}}.
\eeq
Each term of the expression~(\ref{eq:fres}) is responsible for the
scaling~(\ref{eq:csca})
of some of the cumulant amplitudes as $A\to0$.
The first term yields
\beq\label{eq:c0exc}
c_{0,\ell}\approx 2A^2\approx r
\eeq
for all $\ell\ge1$, in agreement with~(\ref{eq:cstir}).
The second term yields
\beq\label{eq:c2exc}
c_{2,0}\to\frac12
\eeq
as $r\to0$, in agreement with~(\ref{eq:c2}).
The third term of~(\ref{eq:fres}) is an odd function of $h$,
and therefore concerns odd values of $k$.
Introducing the series expansion
\beq\label{eq:abdef}
\sqrt{1+\frac{x}{16}}=\sum_{p\ge0}\frac{a_p}{(2p+1)!}\,x^p,\qquad
a_p=(-1)^{p-1}\frac{(2p+1)!}{16^p(2p-1)}\,b_p
\eeq
(see~(\ref{eq:bdef})), i.e.,
\be
a_0=1,\quad
a_1=\frac{3}{16},\quad
a_2=-\frac{15}{256},\quad
a_3=\frac{315}{4096},
\ee
and so on,
we obtain
\beq\label{eq:cres}
c_{2p+1,\ell}\approx a_p\,\left(-\half(2p-1)\right)^\ell\,A^{1-2p}.
\eeq
All cumulant amplitudes $c_{k,\ell}$ with odd $k$ therefore obey the scaling
ansatz~(\ref{eq:csca}),
while the cumulant amplitudes with even $k$ are subleading as $r\to0$,
with the exception of those previously mentioned in~(\ref{eq:c0exc})
and~(\ref{eq:c2exc}).
Finally, the cumulant amplitudes $C_n$ for odd $n$ are governed by the term
$k=n$ in~(\ref{eq:csum}).
We have therefore
\be
C_{2p+1}\approx a_p\,A^{1-2p},
\ee
whereas the $C_n$ with even $n$ are subleading as $r\to0$, except $C_2\approx
c_{2,0}\approx1/2$.

\subsection{Large deviations}

The result~(\ref{eq:s2}) has an alternative interpretation in terms of large
deviations~\cite{ellis,dembo,touchette,vulpiani}.
It implies that the joint probability distribution of $\nc$ and $\nb$
falls off exponentially~as
\beq\label{eq:i2}
\prob(\nc\approx\xi t,\ \nb\approx\eta t)\sim\e^{-I(\xi,\eta)t}
\eeq
in the regime of late times,
for fixed densities $\xi$ of crosses and $\eta$ of dots.
The estimate
\be
\mean{\e^{\lam\nc+\mu\nb}}
\sim\int\dd\xi\int\dd\eta\;\e^{[\lam\xi+\mu\eta-I(\xi,\eta)]t}
\sim\e^{S(\lam,\mu)t}
\ee
shows that the bivariate functions $S(\lam,\mu)$ and $I(\xi,\eta)$
are related by a Legendre transformation of the form
\beq\label{eq:lege}
S(\lam,\mu)+I(\xi,\eta)=\lam\xi+\mu\eta,
\eeq
with
\beq\label{eq:legedifs}
\xi=\frac{\partial S}{\partial\lam},\quad
\eta=\frac{\partial S}{\partial\mu},\quad
\lam=\frac{\partial I}{\partial\xi},\quad
\mu=\frac{\partial I}{\partial\eta}.
\eeq

In the late-time regime,
the joint distribution of $\nc$ and $\nb$ becomes peaked around the point
\be
\xi_0=\ac=A-r=c_{1,0},\qquad\eta_0=\ab=r=c_{0,1},
\ee
in agreement with~(\ref{eq:meannb}),~(\ref{eq:ncasy}) and~(\ref{eq:c1}).
The form of the bivariate large deviation function $I(\xi,\eta)$ around
$(\xi_0,\eta_0)$
is governed by the regime where $\lam$ and $\mu$ are small.
Using the series expansion~(\ref{eq:sers}), we are left with the quadratic
form~(see~(\ref{eq:c2}))
\be
I(\xi,\eta)\approx\frac
{c_{0,2}(\xi-\xi_0)^2-2c_{1,1}(\xi-\xi_0)(\eta-\eta_0)+c_{2,0}(\eta-\eta_0)^2}
{2(c_{2,0}c_{0,2}-c_{1,1}^2)},
\ee
describing the Gaussian bulk of the joint distribution of $\nc$ and $\nb$.

\begin{figure}
\begin{center}
\includegraphics[angle=0,width=0.5\linewidth,clip=true]{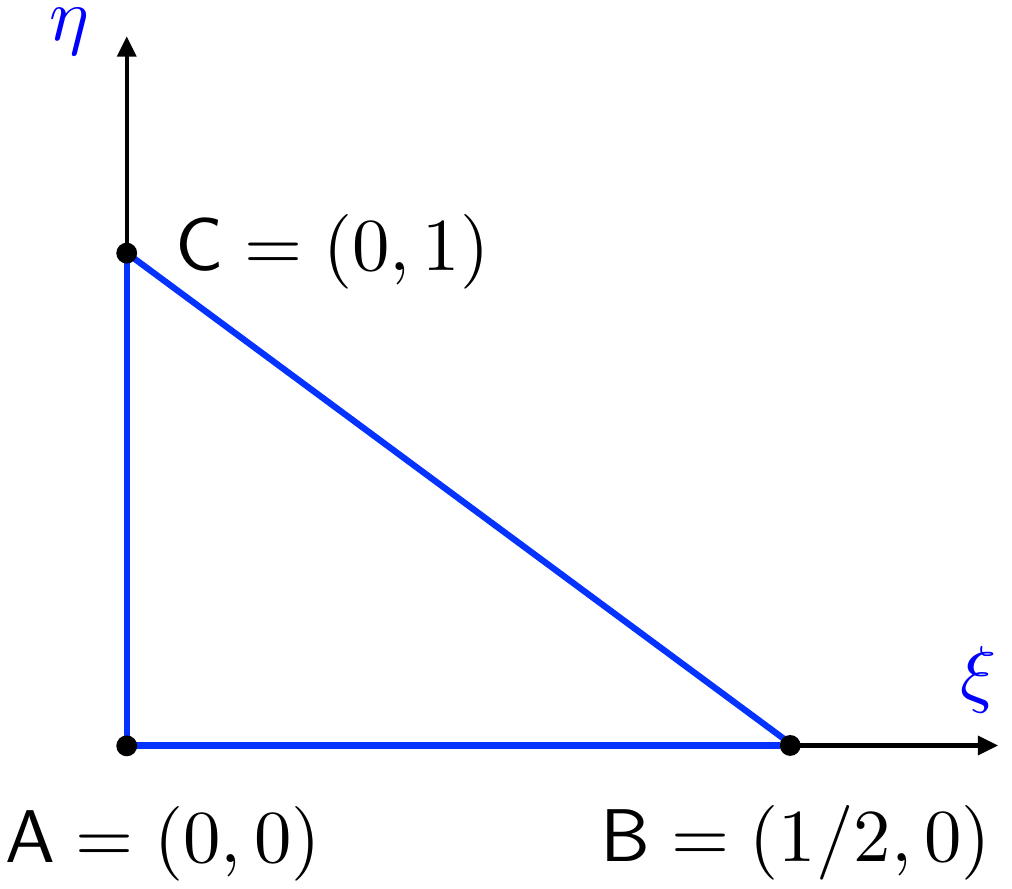}
\caption{\small
Triangular domain of permitted values of densities $\xi$ of returns to the
origin of the walk (crosses) and $\eta$ of resetting events (dots).
}
\label{fig:tri}
\end{center}
\end{figure}

The subsequent analysis shows
that the domain of permitted values of the densities~$\xi$ and $\eta$
is the triangle $\A\B\C$ shown in figure~\ref{fig:tri}.
The large deviation function $I(\xi,\eta)$ is continuous all along the boundary
of the triangle.
Its behaviour near the vertices and the edges of the triangle
is governed by the regime where $\lam$ and/or $\mu$ are large,
either positive or negative.

The vertex $\A=(0,0)$ is reached for $\lam$ and $\mu\to-\infty$.
We obtain $S(-\infty,-\infty)=\ln(1-r)$, and so
\be
I_\A=I(0,0)=-\ln(1-r).
\ee
This can be interpreted as follows.
Point $\A$ corresponds to the situation where the system is empty.
The absence of resettings ($\eta=0$) brings a weight $(1-r)^t$.
Given this condition, the probability of the system containing no crosses
($\xi=0$) is $R(t)$ (see~(\ref{eq:psur})),
which falls off as a power of time (see~(\ref{eq:rasy})),
and thus does not contribute to the large deviation function.
The vertex $\B=(1/2,0)$ is reached for $\lam\to-\infty$ and $\mu\to+\infty$.
We obtain
\be
I_\B=I(1/2,0)=-\ln(1-r)+\frac{1}{2}\,\ln 2.
\ee
The absence of resettings ($\eta=0$) indeed again brings a weight $(1-r)^t$.
Given this condition, when time $t$ is even,
the number $\nc$ of crosses takes its maximal value $t/2$, and $\xi=1/2$,
if the walk consists of a sequence of $t/2$ back-and-forth excursions
on either side of the origin.
This constraint brings a weight $2^{-t/2}$.
The vertex $\C=(0,1)$ is reached for $\lam\to+\infty$ and $\mu\to-\infty$.
We obtain
\be
I_\C=I(0,1)=-\ln r.
\ee
The condition $\eta=1$ indeed amounts to having a resetting event at each time
step.
This brings a weight $r^t$, and there is no space left for crosses.

Along the edge $\A\B$, $I(\xi,0)$ increases monotonically from $I_\A$ to
$I_\B$.
Along $\A\C$, $I(0,\eta)$ is not monotonic and exhibits a minimum,
to be identified with $\ic(0)$ (see below).
Finally, a generic point along $\B\C$ is reached for $\lam$ and $\mu\to+\infty$
in such a way that the difference $\lam-2\mu$ is kept fixed.
The large deviation function thus obtained is not monotonic and exhibits a
minimum.

Let us now turn to the univariate large deviation functions
associated with $\nb$, $\nc$ and their sum $\ncb$.
In the case of $\nb$,~(\ref{eq:i2}) allows us to recover the following results.
We find
\be
\prob(\nb\approx\eta t)\sim\e^{-\ib(\eta)t}\qquad(0<\eta<1),
\ee
where
\be
\ib(\eta)=\min_{\xi}I(\xi,\eta)=\mu\eta-S(0,\mu)
\ee
is the Legendre transform of the function $S(0,\mu)$ given in~(\ref{eq:smu}).
We thus obtain a known expression,
\be
\ib(\eta)=\eta\ln\frac{\eta}{r}+(1-\eta)\ln\frac{1-\eta}{1-r},
\ee
with limit values
\beq\label{{eq:ld1lim}}
\ib(0)=I_\A,\qquad\ib(1)=I_\C,
\eeq
and a quadratic behaviour
\be
\ib(\eta)\approx\frac{(\eta-r)^2}{r(1-r)}
\ee
around $\eta_0=r$.

In the case of $\nc$,~(\ref{eq:i2}) yields
\be
\prob(\nc\approx\xi t)\sim\e^{-\ic(\xi)t}\qquad(0<\xi<1/2),
\ee
where
\be
\ic(\xi)=\min_{\eta}I(\xi,\eta)=\lam\xi-S(\lam,0)
\ee
is the Legendre transform of $S(\lam,0)$.
This function has the limit value
\be
\ic(1/2)=I_\B
\ee
and the quadratic behaviour
\be
\ic(\xi)\approx\frac{(\xi-\ac)^2}{2c_{2,0}}
\ee
around $\xi_0=\ac=c_{1,0}$.
The limit value $\ic(0)$ is given by the decay rate of the distribution of
$\Tz^{(\tt r)}$, introduced in~(\ref{eq:sigrat}),
\be
\ic(0)=\sigma,
\ee
since $\prob(\nc=0)=R^{(\tt r)}(t)\sim \e^{-\sigma t}$, according to~(\ref{eq:dressedR}).

Finally, in the case of $\ncb=\nc+\nb$,~(\ref{eq:i2}) yields
\be
\prob(\ncb\approx\varphi t)\sim\e^{-I(\varphi)t}\qquad(0<\varphi<1),
\ee
where
\be
I(\varphi)=\min_\xi I(\xi,\varphi-\xi)=\lam\varphi-S(\lam,\lam)
\ee
is the Legendre transform of $S(\lam,\lam)$.
The limit values
\be
I(0)=I_\A,\qquad I(1)=I_\C,
\ee
coincide with~(\ref{{eq:ld1lim}}).
The function $I(\varphi)$ has the expected quadratic behaviour
\be
I(\varphi)\approx\frac{(\varphi-A)^2}{2C_2}
\ee
around $\varphi_0=A=C_1$.

\begin{figure}
\begin{center}
\includegraphics[angle=0,width=0.7\linewidth,clip=true]{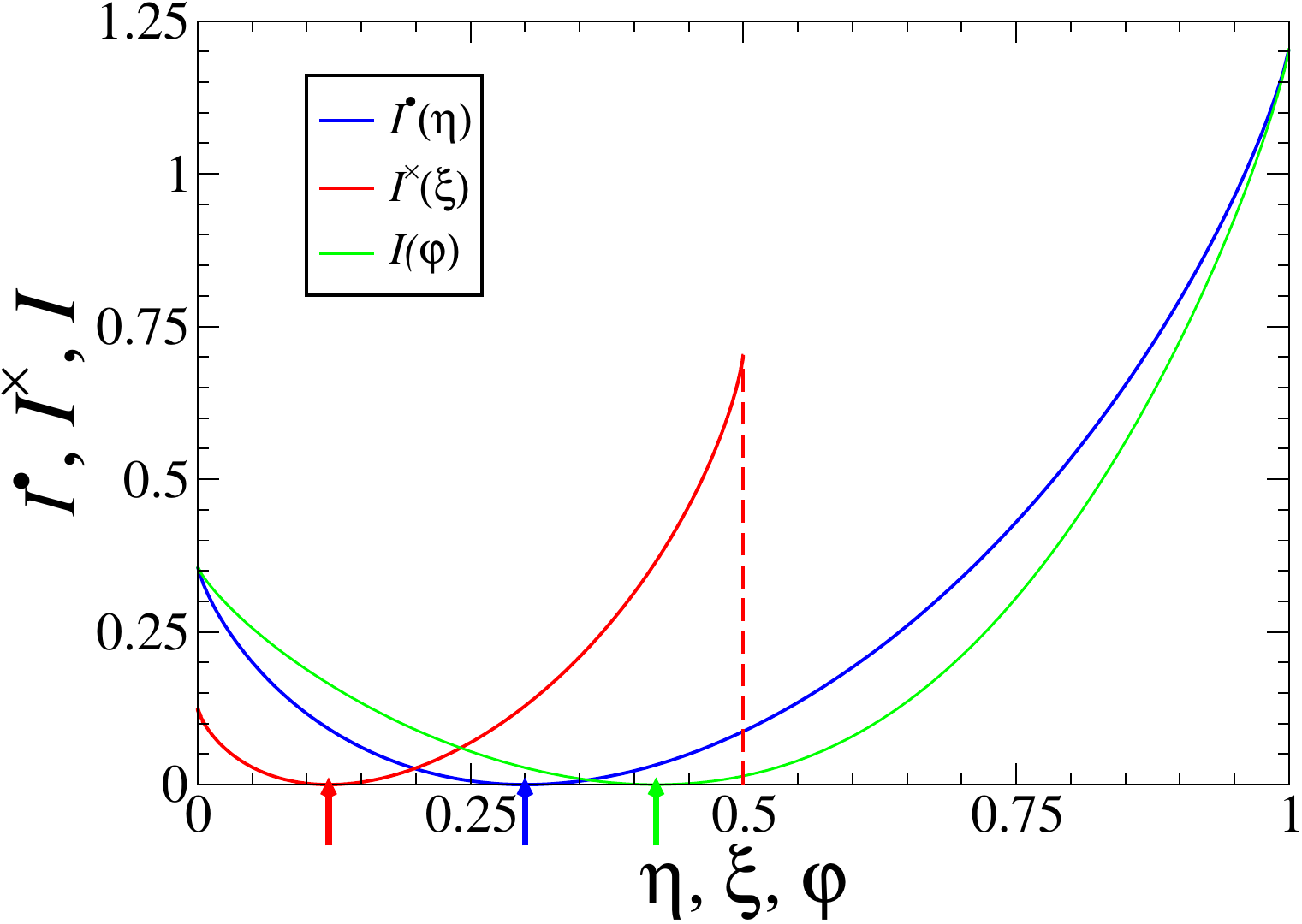}
\caption{\small
Univariate large deviation functions $\ib(\eta)$, $\ic(\xi)$, and $I(\varphi)$,
corresponding to $\nb$, $\nc$ and their sum $\ncb$, respectively,
for a resetting probability $r=0.3$,
so that $\eta_0=0.3$, $\xi_0=0.120084\dots$ and $\varphi_0=0.420084\dots$
(arrows).}
\label{fig:lds}
\end{center}
\end{figure}

Figure~\ref{fig:lds} shows plots of the univariate large deviation functions
$\ib(\eta)$, $\ic(\xi)$, and $I(\varphi)$,
respectively corresponding to $\nb$, $\nc$ and their sum $\ncb$,
for a resetting probability $r=0.3$.
Numerical values of these functions are obtained by means of~(\ref{eq:s2trig}).
All derivatives required by the Legendre
transform~(\ref{eq:lege}),~(\ref{eq:legedifs})
are worked out by analytical means.
For instance,
\be
\xi
=\frac{\partial S}{\partial\lam}
=-\frac{\partial\ln w_1}{\partial\ln z}
=-\frac{z}{w_1}\frac{\partial w_1}{\partial z}
=\left.\frac{z\,\partial P/\partial z}{w\,\partial P/\partial
w}\right\vert_{w=w_1}
\ee
(see~(\ref{eq:wpol})) is a rational expression in $r$, $z$, $y$ and $w_1$.

\section{Crossover regime at weak resetting}
\label{sec:cross}

The statistics of the number $\nc$ of crosses exhibits a non-trivial behaviour
in the crossover regime of weak resetting ($r\to0$) and late times
($t\to\infty$).
In the absence of resettings, the mean value $\mean{\nc}$ scales as $\sqrt{t}$
(see~(\ref{eq:Ntasym})), while in the case of weak resetting, it scales as
$\sqrt{r}\,t$ (see~(\ref{eq:ncasy}) and~(\ref{eq:acres})).
These two estimates become comparable when the product $rt$ is of order unity.
Interestingly, the latter is precisely the value of the mean number of
resettings $\mean{\nb}$ (see~(\ref{eq:meannb})), which implies that a finite
number of resetting events are sufficient to induce a macroscopic crossover in
the statistics of $\nc$.
This phenomenon has also been described in other observables, including the
maximum and number of records of random walks under weak
resetting~\cite{m2s2,glmax}.

The full distribution of $\nc$ throughout this crossover regime
can be derived from~(\ref{eq:key}).
Setting $w=\e^{-s}$, $y=1$ and $z=\e^{-p}$,
and working to leading order in the continuum regime where $r$, $s$ and $p$ are
small,
we obtain
\be
\int_0^\infty\dd
t\,\e^{-st}\mean{\e^{-p\nc}}\approx\frac{1}{s+p\sqrt{(r+s)/2}}.
\ee
Inverting the Laplace transform in $p$ yields
\beq\label{{eq:rlapn}}
\int_0^\infty\dd t\,\e^{-st}\,\prob(\nc=\nn)
\approx\sqrt\frac{2}{r+s}
\exp\Big(-s\sqrt\frac{2}{r+s}\,\nn\Big).
\eeq
The two expressions above are very similar to~(\ref{eq:laplapn})
and~(\ref{eq:lapn}).
We thus infer from~(\ref{{eq:rlapn}}) that the number $\nc$ of crosses scales
as
\beq\label{{eq:xidef}}
\nc\approx\sqrt{t}\,\mzeta,
\eeq
where the rescaled random variable $\mzeta$ has a limiting distribution with
density $f(\zeta,u)$,
depending only on the parameter
\be
u=rt=\mean{\nb}.
\ee
Introducing the ratio $\lam=s/r$,~(\ref{{eq:rlapn}}) becomes
\beq\label{eq:phidef}
f(\zeta,u)=\sqrt{2u}\int\frac{\dd\lam}{2\pi\ii}\frac{\e^{\lam
u}}{\sqrt{\lam+1}}
\exp\left(-\frac{\lam}{\sqrt{\lam+1}}\,\sqrt{2u}\,\zeta\right).
\eeq
This expression cannot be made more explicit,
except at $\zeta=0$, resulting in the following value:
\be
f(0,u)=\sqrt{2u}\int\frac{\dd\lam}{2\pi\ii}\frac{\e^{\lam u}}{\sqrt{\lam+1}}
=\sqrt\frac{2}{\pi}\,\e^{-u}.
\ee
Hereafter we examine the behaviour of this distribution
in the regimes where the parameter $u$ is either small or large.
We then shift our focus to the analysis of the moments and the cumulants of
$\mzeta$.

\subsubsection*{Behaviour for $u\ll1$}

The behaviour of $f(\zeta,u)$ for small $u$ can be derived
by setting $\lambda=p^2/u$ in~(\ref{eq:phidef}),
and expanding the integrand as a power series in~$u$ at fixed $p$.
We thus obtain
\beqa
f(\zeta,u)
&=&\int\frac{\dd p}{2\pi\ii}\,\e^{p^2-\sqrt{2}p\zeta}
\left[2\sqrt{2}+\left(\frac{2\zeta}{p}-\frac{\sqrt{2}}{p^2}\right)u+\cdots\right]
\nonumber\\
&=&\sqrt\frac{2}{\pi}\,\e^{-\zeta^2/2}
+\left(2\zeta\erfc\frac{\zeta}{\sqrt{2}}-\sqrt\frac{2}{\pi}\,\e^{-\zeta^2/2}\right)u+\cdots,
\label{{eq:usmall}}
\eeqa
where erfc is the complementary error function.
The first term matches the half-Gaussian asymptotic
distribution~(\ref{eq:hgau}) of $\nc$ in the absence of resetting.

\subsubsection*{Behaviour for $u\gg1$}

Let us define the random variable $X$ by
\beq\label{eq:xsplit}
\mzeta=\sqrt\frac{u}{2}+X.
\eeq
The behaviour of $f(\zeta,u)$ for large $u$ can then be derived by setting
$\lambda=p/\sqrt{2u}$ in~(\ref{eq:phidef})
and expanding the integrand as a power series in~$1/\sqrt{u}$ at fixed $p$.
We thus obtain, with $\zeta=\sqrt{u/2}+x$,
\beqa
f(\zeta,u)
&=&\int\frac{\dd p}{2\pi\ii}\,\e^{p^2/4-px}
\left(1+\frac{8p^2x-3p^3-8p}{16\sqrt{2u}}+\cdots\right)
\nonumber\\
&=&\frac{\e^{-x^2}}{\sqrt\pi}
\left(1+\frac{x(2x^2+1)}{4\sqrt{2u}}+\cdots\right).
\label{{eq:ularge}}
\eeqa
Whenever the parameter $u$ is large,
the random variable $\mzeta$ therefore consists of a large deterministic part,
growing as $\sqrt{u}$,
along with a fluctuating component $X$ of order unity.
To leading order, the distribution of $X$ is a Gaussian with variance $1/2$.

\begin{figure}
\begin{center}
\includegraphics[angle=0,width=0.7\linewidth,clip=true]{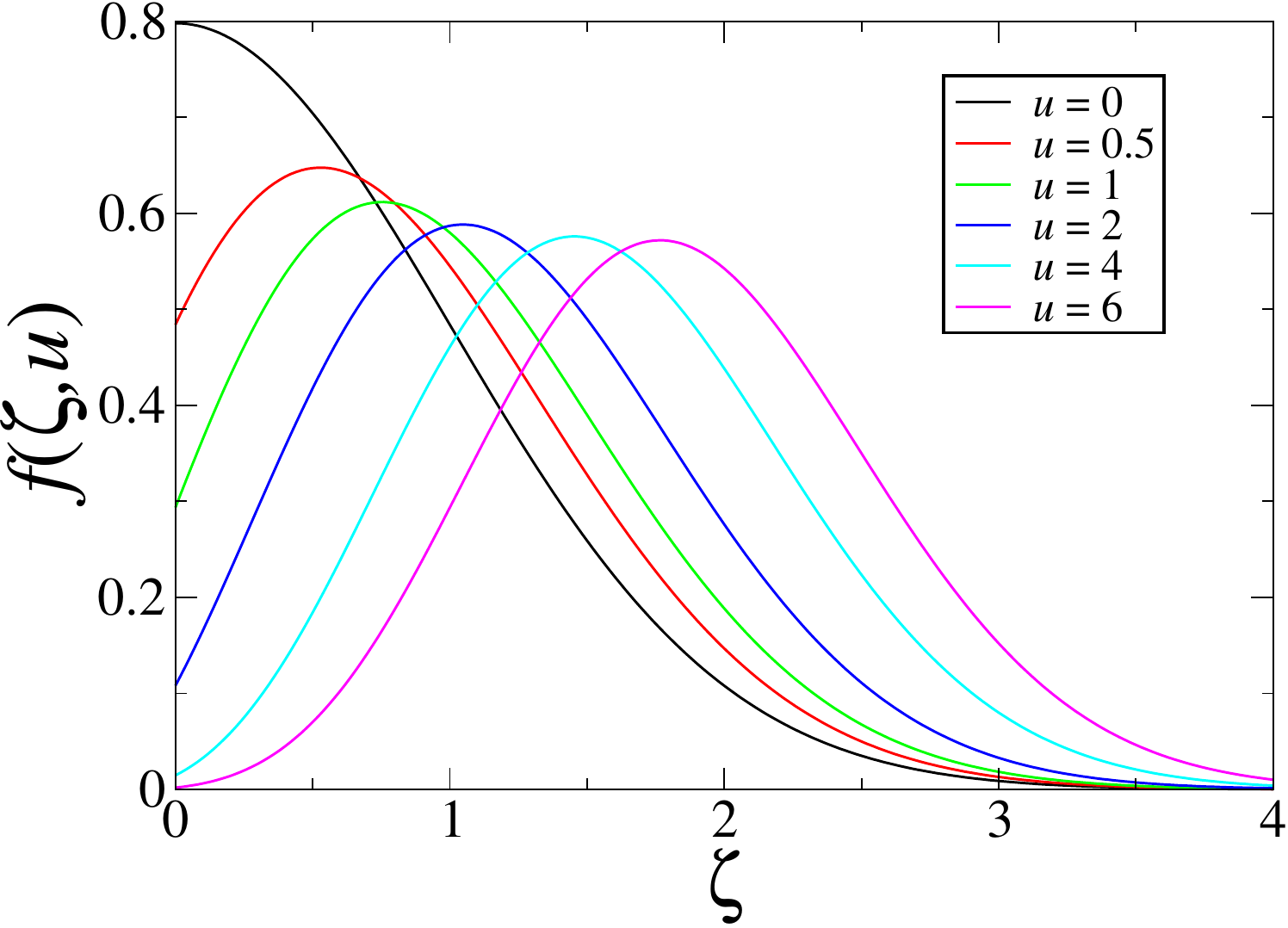}
\caption{\small
Distribution $f(\zeta,u)$ of the rescaled number $\zeta$ of returns to the
origin in the weak-resetting crossover regime (see~(\ref{{eq:xidef}})),
for several values of the parameter $u=rt$ (see legend).}
\label{fig:fzeta}
\end{center}
\end{figure}

Figure~\ref{fig:fzeta} shows the distribution $f(\zeta,u)$
for several values of the parameter $u=rt$ (see legend).
The plotted distributions exhibit a smooth but rather rapid crossover
between a half-Gaussian form at $u=0$ (see~(\ref{{eq:usmall}}))
and a shifted Gaussian at large $u$ (see~(\ref{{eq:ularge}})).
In particular, the maxima of the curves converge very fast to their limit
$1/\sqrt{\pi}=0.564189\dots$

\subsubsection*{Moments of $\mzeta$}

Equation~(\ref{eq:phidef}) yields the following formula for the moments
of~$\zeta$:
\beq\label{eq:momint}
\mu_k(u)=\mean{\mzeta^k}
=\int_0^\infty\dd\zeta\, \zeta^k\,f(\zeta,u)
=\frac{k!}{(2u)^{k/2}}\int\frac{\dd\lam}{2\pi\ii}\,\e^{\lam
u}\,\frac{(\lam+1)^{k/2}}{\lam^{k+1}}.
\eeq
These moments only depend on the parameter $u=rt$.
They are such that
\be
\mean{(\nc)^k}\approx\mu_k(u)\,t^{k/2}
\ee
throughout the crossover regime.
The moments~(\ref{eq:hgaumoms}) in the absence of resetting yield
\beq\label{eq:momszero}
\mu_{2n}(0)=\frac{(2n)!}{2^nn!},\qquad
\mu_{2n+1}(0)=\sqrt\frac{2}{\pi}\,2^nn!.
\eeq
The results above suggest that the moments $\mu_k(u)$
have different analytical expressions according to the parity of the integer
exponent $k$.
This is indeed the case.

For even $k=2n$,
the integrand in the rightmost side of~(\ref{eq:momint}) is a rational function
of $\lam$.
Expanding out $(\lam+1)^n$, we readily obtain
\be
\mu_{2n}(u)=\frac{(2n)!n!}{2^n}\sum_{m=0}^n\frac{u^m}{m!(n-m)!(n+m)!},
\ee
namely
\be
\mu_2(u)=\frac{u+2}{2},\quad
\mu_4(u)=\frac{u^2+8u+12}{4},\quad
\mu_6(u)=\frac{u^3+18u^2+90u+120}{8},
\ee
and so on.
The moment $\mu_{2n}(u)$ is a polynomial of degree $n$ in $u$.
Its constant term $(m=0)$ matches the first expression in~(\ref{eq:momszero}),
whereas its leading term $(m=n)$ yields
\be
\mu_{2n}(u)\approx\left(\frac{u}{2}\right)^n,
\ee
in agreement with~(\ref{eq:xsplit}).

For odd $k=2n+1$,
the integrand in the rightmost side of~(\ref{eq:momint}) is now the ratio
of a rational function by $\sqrt{\lam+1}$.
Proceeding as before, we obtain
\be
\mu_{2n+1}(u)=\frac{(2n+1)!(n+1)!}{(2u)^{n+1/2}}
\sum_{m=0}^{n+1}\frac{g_{n+m}(u)}{m!(n+1-m)!(n+m)!},
\ee
with
\be
g_n(u)=n!\int\frac{\dd\lam}{2\pi\ii}\,
\frac{\e^{\lam u}}{\lam^{n+1}\sqrt{\lam+1}}
=\int_0^u \dd v\,(u-v)^n\,\frac{\e^{-v}}{\sqrt{\pi v}}.
\ee
It can be shown using two integrations by parts
that these functions obey the three-term linear recursion
\be
g_n(u)=(u-n+\half)g_{n-1}(u)+(n-1)ug_{n-2}(u),
\ee
with initial values
\be
g_0(u)=\erf\sqrt{u},\qquad
g_1(u)=(u-\half)\erf\sqrt{u}+\sqrt\frac{u}{\pi}\,\e^{-u}.
\ee
We thus obtain
\bea
\mu_1(u)&=&(2u+1)\frac{\erf\sqrt{u}}{2\sqrt{2u}}+\frac{\e^{-u}}{\sqrt{2\pi}},
\nonumber\\
\mu_3(u)&=&(8u^3+36u^2+18u-3)\frac{\erf\sqrt{u}}{16u\sqrt{2u}}
+(4u^2+16u+3)\frac{\e^{-u}}{8u\sqrt{2\pi}},
\eea
and so on.
The general structure of the odd moments emerges from the above examples.
Their values at $u=0$ (see~(\ref{eq:momszero})) cannot be easily read off,
as more and more compensations are involved in taking the $u\to0$ limit.
To leading order at large $u$, we have
\be
\mu_{2n+1}(u)\approx\left(\frac{u}{2}\right)^{n+1/2},
\ee
again in agreement with~(\ref{eq:xsplit}).

To close, we mention that the probability $U^{(\tt r)}(t)$ introduced in~(\ref{eq:udef})
scales as
\beq
U^{(\tt r)}(t)\approx\frac{G(u)}{\sqrt{t}}
\label{usca}
\eeq
throughout the crossover regime,
with
\beq
G(u)=\half\mu_1(u)+u\mu_1'(u)
=\frac{\sqrt{u}\,\erf\sqrt{u}}{\sqrt{2}}+\frac{\e^{-u}}{\sqrt{2\pi}}.
\eeq
The probability $U^{(\tt r)}(t)$ exhibits even-odd oscillations,
so that~(\ref{usca}) actually describes the behaviour of the local average
$\half(U^{(\tt r)}(t)+U^{(\tt r)}(t-1))$.

\subsubsection*{Cumulants of $\mzeta$}

In order to compare the above analysis of the crossover with the outcomes of
section~\ref{sec:deviations},
let us consider the cumulants
\be
\gamma_k(u)=\mean{\mzeta^k}_c.
\ee
At large values of $u$,
neglecting exponentially small corrections, these quantities read
\bea
\gamma_1(u)&\approx&\frac{2u+1}{2\sqrt{2u}},\qquad
\gamma_2(u)\approx\frac{4u-1}{8u},
\nonumber\\
\gamma_3(u)&\approx&\frac{6u-1}{16u\sqrt{2u}},\qquad
\gamma_4(u)\approx\frac{3}{32u^2},
\nonumber\\
\gamma_5(u)&\approx&-\frac{3(10u-3)}{128u^2\sqrt{2u}},\qquad
\gamma_6(u)\approx-\frac{15}{64u^3}.
\eea
The cumulants of $\mzeta$ appear to have a simpler dependence on $u$
than the corresponding moments.
To leading order as $u\gg1$, the odd cumulants scale as
\be
\gamma_{2n+1}(u)\approx a_n\left(\frac{u}{2}\right)^{1/2-n},
\ee
in agreement with~(\ref{eq:cres}),
where the amplitudes $a_n$ are given in~(\ref{eq:abdef}).
The second cumulant (variance) admits a finite limit $1/2$,
to be identified with the limit of $c_{2,0}$,
whereas higher even cumulants scale as
\be
\gamma_{2n}(u)\approx\frac{\alpha_n}{u^n},
\ee
for some constants $\alpha_n$.
They are therefore subleading with respect to~the odd ones.

\section{Discussion}
\label{sec:discuss}

To conclude, let us put the main outcomes of the present work
in perspective with those of the companion paper~\cite{nested}.
The point process considered in this latter work involves two generic nested
renewal processes,
an internal one characterised by the distribution $\rho(\tau)$ of interarrival
times,
and an external one characterised by the distribution $f(T)$ of time intervals
between resetting events.
In~\cite{nested}, the main emphasis was on the number $\nc$
of (internal) renewal events occurring within a fixed observation time $t$.
The statistics of this observable revealed a wide variety of asymptotic
behaviours,
dependent on the values of the exponents $\theta_1$ and $\theta_2$
governing the tails of the distributions $\rho(\tau)$ and $f(T)$.
These findings highlight the dominance of the more regular of the two
processes,
specifically the one with the larger tail exponent,
$\tilde\theta=\max(\theta_1,\theta_2)$.
More specifically,
$\nc$ grows linearly in time and has relatively negligible fluctuations
whenever $\tilde\theta>1$,
whereas $\nc\sim t^{\tilde\theta}$ grows subextensively over time while
continuing to fluctuate
for $\tilde\theta<1$.

The reset P\'olya walk considered in the present work is a specific instance of
the general process made of two arbitrary nested renewal processes.
The internal renewal process describes
the spontaneous returns of the walker to its starting point,
whereas the external one consists of resettings, taking place with probability
$r$ at each time step.
In the phase diagram of~\cite{nested},
this example corresponds to $\theta_1=1/2$ and $\theta_2=\infty$,
and hence $\tilde\theta=\infty$,
so that a high degree of regularity is expected for the entire process.

The present analysis corroborates this prediction and completes it by a breadth
of quantitative results
concerning the joint statistics of the numbers $\nc$ of crosses (spontaneous
returns)
and $\nb$ of dots (resetting events) in the regime of late times.
The most salient of these outcomes---highlighted in the
introduction---concerning the linear growth of all joint cumulants
$\mean{(\nc)^k(\nb)^\ell}_c$
and the smoothness of the bivariate large deviation function $I(\xi,\eta)$,
testify that the numbers $\nc$ and $\nb$ are extensive in a strong sense,
and that the reset P\'olya walk indeed manifests a very high degree of regularity.
This characteristic can be related to the exponentially decaying, hence strongly localised,
steady-state distribution of the walker's position under stochastic resetting
(see~(\ref{pst})).
It would be worth investigating whether similar regularity properties
also manifest themselves in other observables pertaining to the P\'olya walk
under resetting.

\subsubsection*{Data availability statement}

The authors have no data to share.

\subsubsection*{Conflict of interest}

The authors declare no conflicts of interest.

\begin{appendices}

\section{Stationary distribution of the walker's position}
\label{sec:appb}

In this appendix we give a brief self-contained derivation of the expression~(\ref{pst})
of the distribution $p(x)$ of the walker's position $x$
in the nonequilibrium steady state
reached in the limit of infinitely large times in the presence of resetting.
The latter expression had been derived in a general setting in~\cite{glmax}.

The stochastic recursion~(\ref{eq:defS}) obeyed by the walker's position
translates into the following balance equation for the corresponding stationary distribution $p(x)$:
\beq
p(x)=r\delta(x,0)+\frac{1-r}{2}(p(x+1)+p(x-1)).
\label{eqpx}
\eeq
The first step in solving this equation consists in
looking for an exponential solution of the form $\mu^x$ away from the origin.
This yields the following quadratic equation for $\mu$:
\beq
(1-r)\mu^2-2\mu+1-r=0,
\eeq
whose solutions read $\mu=\lam$ and $\mu=1/\lam$,
where the larger solution $\lam>1$ reads
\beq
\lam=\frac{1+\sqrt{r(2-r)}}{1-r}.
\eeq
The second step consists in looking for a normalized global solution to~(\ref{eqpx}).
We thus find
\beq
p(x)=\sqrt\frac{r}{2-r}\,\lam^{-\abs{x}},
\eeq
that is~(\ref{pst}).

\section{Zeros of the denominator of~(\ref{eq:key})}
\label{sec:appa}

This appendix is devoted to the determination of the zeros of the denominator
of~(\ref{eq:key})
in the variable~$w$.
Using~(\ref{eq:Zzw}) and~(\ref{eq:fctngen}),
and eliminating the square root entering the formula thus obtained,
results in a polynomial equation of degree three for the zeros $w(z,y)$,
reading
\beq\label{eq:wpol}
P(r,z,y,w)=P_3w^3+P_2w^2+P_1w+P_0=0,
\eeq
with
\beqa
P_3&=&(1-r)((1-r)z+ry)^2,
\qquad
P_2=r^2y^2-(1-r)^2z^2,
\nonumber\\
P_1&=&(1-r)(1-2z)-2r yz,
\qquad
P_0=2z-1.
\label{eq:pdefs}
\eeqa
Polynomial equations of degree three can be solved analytically,
either by Cardano's method or by the trigonometric method.
Here, we adopt the latter approach,
which has the advantage that no complex numbers are involved
when the three zeros are real.
This is indeed the case here, for small enough real $\lam$ and $\mu$.
Setting $w=B+x$, we arrive at a reduced equation for $x$,
\be
x^3+px+q=0,
\ee
with
\beq\label{eq:bpq}
B=-\frac{P_2}{3P_3},
\quad
p=\frac{P_1}{P_3}-\frac{P_2^2}{3P_3^2},
\quad
q=\frac{P_0}{P_3}-\frac{P_1P_2}{3P_3^2}+\frac{2P_2^3}{27P_3^3}.
\eeq
The condition for all zeros to be real reads $4p^3+27q^2<0$,
implying in particular $p<0$.
These zeros then read
\beqa
w_1&=&B+\sigma\cos(\theta-2\pi/3),
\nonumber\\
w_2&=&B+\sigma\cos\theta,
\nonumber\\
w_3&=&B+\sigma\cos(\theta+2\pi/3),
\label{eq:wres}
\eeqa
with
\beq\label{eq:sigth}
\sigma=-\sqrt{-\frac{4p}{3}},\qquad
\cos3\theta=-\frac{4q}{\sigma^3}\qquad
(0\le\theta\le\pi/3).
\eeq
For small enough real $\lam$ and $\mu$,
the smallest of the three zeros is $w_1$, which is positive, so that
\be
S(\lam,\mu)=-\ln w_1,
\ee
that is~(\ref{eq:s2trig}).

\end{appendices}

\bibliography{paperPolya.bib}

\end{document}